\documentclass[10pt]{amsart}
\usepackage{graphicx,amssymb}
\usepackage{a4wide}
\usepackage{xypic}
\xyoption{all}
\xyoption{arc}
\xyoption{poly}

\theoremstyle{plain}
\newtheorem{theorem}{Theorem}[section]
\newtheorem{lemma}[theorem]{Lemma}
\newtheorem{corollary}[theorem]{Corollary}
\newtheorem{proposition}[theorem]{Proposition}

\theoremstyle{definition}
\newtheorem{definition}[theorem]{Definition}
\newtheorem{example}[theorem]{Example}

\newtheorem{remark}[theorem]{Remark}

\newcommand{\thickhline}{\noalign{\hrule height 0.9pt}}

\def\arA{\mathbf A}
\def\arB{\mathbf B}
\def\arC{\mathbf C}
\def\arD{\mathbf D}
\def\arI{\mathbf I}

\def\Z{\mathbb Z}
\def\Q{\mathbb Q}

\def\c{\mathbf c}
\def\d{\mathbf d}

\def\t{\operatorname{\it t}}
\def\len{\operatorname{len}}
\def\INF{\t_{\inf}}
\def\SUP{\t_{\sup}}
\def\LEN{\t_{\len}}
\def\infs{\inf\nolimits{\!}_s}
\def\sups{\sup\nolimits{\!}_s}
\def\lens{\len{\!}_s}

\def\wedgeL{\wedge_L}
\def\wedgeR{\wedge_R}

\def\myangle#1{\langle #1\rangle}

\begin{document}

\title{Periodic elements in Garside groups}

\author{Eon-Kyung Lee and Sang-Jin Lee}
\address{Department of Mathematics, Sejong University,
    Seoul, 143-747, Korea}
\email{eonkyung@sejong.ac.kr}
\address{Department of Mathematics, Konkuk University,
    Seoul, 143-701, Korea}
\email{sangjin@konkuk.ac.kr}

\begin{abstract}
Let $G$ be a Garside group with Garside element $\Delta$,
and let $\Delta^m$ be the minimal positive central power of $\Delta$.
An element $g\in G$ is said to be \emph{periodic} if
some power of it is a power of $\Delta$.
In this paper, we study periodic elements in Garside groups
and their conjugacy classes.

We show that
the periodicity of an element does not depend on
the choice of a particular Garside structure
if and only if the center of $G$ is cyclic;
if $g^k=\Delta^{ka}$ for some nonzero integer $k$,
then $g$ is conjugate to $\Delta^a$;
every finite subgroup of the quotient group
$G/\langle \Delta^m\rangle$ is cyclic.

By a classical theorem of Brouwer, Ker\'ekj\'art\'o and Eilenberg,
an $n$-braid is periodic if and only if
it is conjugate to a power of one of two specific
roots of $\Delta^2$.
We generalize this to Garside groups by showing
that every periodic element is conjugate to a power of a root of $\Delta^m$.

We introduce the notions of slimness and precentrality for periodic elements,
and show that the super summit set of a slim, precentral periodic element
is closed under any partial cycling.
For the conjugacy problem, we may assume the slimness without loss of generality.
For the Artin groups of type $\arA_n$, $\arB_n$, $\arD_n$, $\arI_2(e)$ and
the braid group of the complex reflection group of type $(e,e,n)$,
endowed with the dual Garside structure, we may further assume the precentrality.

\end{abstract}

\maketitle


\section{Introduction}

Garside groups, first introduced by Dehornoy and Paris~\cite{DP99},
are a lattice-theoretic generalization of braid groups
and Artin groups of finite type.
Garside groups are equipped with a special element $\Delta$,
called the \emph{Garside element}.
An element $g$ of a Garside group is said to be \emph{periodic} if
$$
g^k=\Delta^\ell
$$
for some integers $k\neq 0$ and $\ell$~\cite{Bes06b,BGG08}.

Recently there were several results on periodic elements of Garside groups
such as the characterization of finite subgroups
of the central quotient of finite type Artin groups by Bestvina~\cite{Bes99}
and its extension to Garside groups by Charney, Meier and Whittlesey~\cite{CMW04};
the characterization of periodic elements in the braid groups of
complex reflection groups by Bessis~\cite{Bes06a};
a new algorithm for solving the conjugacy search problem for periodic braids
by Birman, Gebhardt and Gonz\'alez-Meneses~\cite{BGG07}.

In this paper we study periodic elements in Garside groups.
We are interested in general Garside groups,
but also concerned with particular Garside groups such as
the Artin groups $A(\arA_n)$, $A(\arB_n)$, $A(\arD_n)$
and $A(\arI_2(e))$,
and the braid group $B(e,e,n)$ of
the complex reflection group of type $(e,e,n)$.

\subsection{Periodicity and Garside structure}
The periodicity of an element in a Garside group
generally depends on the choice of a particular Garside structure,
more precisely on the Garside element.
A Garside group may admit more than one Garside structure.
Therefore it is natural to ask the following question.

\begin{quote}
When is the periodicity of an element
independent of the choice of a particular Garside structure?
\end{quote}

It is easy to see that the periodicity does not depend on the choice of a Garside structure
if and only if any two Garside elements are commensurable,
and that this happens if the center is cyclic.
(Two elements $g$ and $h$ of a group are said to be {\em commensurable}
if $g^k$ is conjugate to $h^\ell$ for some nonzero integers $k$ and $\ell$.)
We show that the converse is also true.

\medskip\noindent
\textbf{Theorem~\ref{thm:Gar_struc}.}\ \ \em
Let $G$ be a Garside group.
Then the center of $G$ is cyclic if and only if
any pair of Garside elements of $G$ are commensurable.
\upshape
\medskip

The irreducible Artin groups of finite type and, more generally,
the braid groups of irreducible well-generated complex reflection groups
are Garside groups with cyclic center~\cite{Bes06a}.
Therefore, in these groups, an element is periodic (with respect to a Garside element)
if and only if it has a central power.
However, not all Garside groups have cyclic center.
A typical example is $\Z^\ell$ for $\ell\ge 2$ (see Example~\ref{ex:Per_GarStru}).

\subsection{Roots of periodic elements}\label{ssec:Intr_Root}

We begin with a definition of Bessis in~\cite{Bes06b}:
for a Garside group $G$ with Garside element $\Delta$,
an element $g\in G$ is $p/q$-periodic if $g^q=\Delta^p$ for $p\in\Z$ and $q\in\Z_{\ge 1}$.

Note that
$g^k=\Delta^\ell$ for some $k\in\Z_{\ge 1}$ and $\ell\in\Z$
if and only if
$g^q$ is conjugate to $\Delta^p$ for some $q\in\Z_{\ge 1}$ and $p\in\Z$,
because $\Delta$ has a central power.
Using this equivalence, we define the notion of $p/q$-periodicity
in a slightly different way.

\medskip\noindent
\textbf{Definition~\ref{defn:pq-per}.}\ \
Let $G$ be a Garside group with Garside element $\Delta$.
An element $g\in G$ is said to be \emph{$p/q$-periodic}
for $p\in\Z$ and $q\in\Z_{\ge 1}$ if $g^q$ is conjugate to $\Delta^p$ and
$q$ is the smallest positive integer such that $g^q$
is conjugate to a power of $\Delta$.
\medskip

In the above definition, the $p/q$-periodicity \emph{a priori}
depends on the actual $p$ and $q$
and not just on the rational number $p/q$ because it may happen that
$g^{kq}$ is conjugate to $\Delta^{kp}$ for some $k\ge 2$
but $g^q$ is not conjugate to $\Delta^p$.
Motivated by this observation, we show the following.

\medskip\noindent
\textbf{Theorem~\ref{thm:unique}.}\ \ \em
Let $G$ be a Garside group with Garside element $\Delta$, and let $g\in G$
and $a, b, k \in\Z_{\neq 0}$.
\begin{itemize}
\item[(i)]
If\/ $g^{kb}$ is conjugate to $\Delta^{ka}$,
then $g^b$ is conjugate to $\Delta^a$.

\item[(ii)]
If\/ each of $g^{a}$ and $g^{b}$ is conjugate to a power of $\Delta$,
then so is $g^{a\wedge b}$, where $a\wedge b$ denotes
the greatest common divisor of $a$ and $b$.
\end{itemize}
\upshape
\medskip

By this theorem, the term `$p/q$-periodic' contains that
$p$ and $q$ are coprime.

The above theorem is a sort of uniqueness property of roots up to conjugacy.
On this property, stronger results are known for some specific groups.
Let $g$ and $h$ be elements of a group $G$ with
\begin{equation}\label{eq:root}
g^k=h^k\qquad\mbox{for some $k\ne 0$}.
\end{equation}
If $G$ is the pure $n$-braid group $P_n$, then
$g=h$ by Bardakov~\cite{Bar92}.
(This also follows from the biorderability of the pure braid groups
by Kim and Rolfsen~\cite{KR03}.)
If $G$ is the $n$-braid group $B_n$, then $g$ and $h$ are conjugate
by Gonz\'alez-Meneses~\cite{Gon03}.
If $G$ is the Artin group of type $\arB$, $\tilde \arA$ or $\tilde \arC$,
then $g$ and $h$ are conjugate~\cite{LL10}.
If $G$ is the braid group of a well-generated complex reflection group
and $g$ and $h$ are periodic elements,
then $g$ and $h$ are conjugate by Bessis~\cite{Bes06a}.
For a study of roots in mapping class groups, see~\cite{BP09}.

Theorem~\ref{thm:unique} shows that if $G$ is a Garside group
and $h$ is a power of a Garside element $\Delta$, then (\ref{eq:root}) implies that
$g$ and $h$ are conjugate.
In Garside groups, even for periodic elements, it is hard to obtain a result
stronger than Theorem~\ref{thm:unique}.
For every $k\ge 2$, there is a Garside group with periodic elements $g$ and $h$
such that $g^k=h^k$ but $g$ and $h$ are not conjugate.
(See Example~\ref{ex:nonunique}.)

\medskip
The following is a question of Bessis~\cite[Question~4]{Bes06b}.

\medskip
\noindent\textbf{Question.}\ \
Let $G$ be a Garside group with Garside element $\Delta$.
Let $g\in G$ be a periodic element with respect to $\Delta$.
Does $G$ admit a Garside structure with Garside element $g$?
\medskip

The above question is answered almost positively in the case of the braid group $B_n$:
each periodic element in $B_n$ is conjugate to
a power of one of the particular braids $\delta$ and $\varepsilon$
which are the Garside elements in the dual Garside structures of $B_n$
and $A(\mathbf B_{n-1})$, respectively, where $A(\mathbf B_{n-1})$ denotes
the Artin group of type $\mathbf B_{n-1}$ viewed as a subgroup of $B_n$.
In~\cite{Bes06b}, Bessis showed that the above question is answered almost positively
in the setting of Garside groupoids.

To a Garside group $G$ with an affirmative answer to the above question,
the idea of Birman, Gebhardt and Gonz\'alez-Meneses in~\cite{BGG07}
can possibly be applied.
Precisely, in order to solve the conjugacy search problem for periodic elements
$g$ and $h$ of $G$,
it suffices to find a Garside structure with Garside element $g$.

Using Theorem~\ref{thm:unique}, we give a negative answer to the above question:
there is a Garside group $G$ with a periodic element $g$
such that there is no Garside structure on $G$ with Garside element $g$.
(See Example~\ref{ex:nonunique}.)

\subsection{Finite subgroups of the quotient group $G_\Delta$}

In a Garside group $G$, the Garside element $\Delta$ always has a central power.
Let $\Delta^m$ be the minimal positive central power.
Let $G_\Delta$ be the quotient $G/\langle\Delta^m\rangle$,
where $\myangle{\Delta^m}$ is the cyclic group generated by $\Delta^m$.
For an element $g\in G$, let $\bar g$ denote the image of $g$
under the natural projection from $G$ to $G_\Delta$.
Hence, an element $g\in G$ is periodic if and only if
$\bar g$ has a finite order in $G_\Delta$.

About finite subgroups of $G_\Delta$, the following facts are known.

\begin{itemize}
\item[(i)]
If $G$ is an Artin group of finite type,
then every finite subgroup of\/ $G_\Delta$ is cyclic.

\item[(ii)]
If $G$ is a Garside group, then every finite subgroup of\/ $G_\Delta$
is abelian of rank at most 2.
\end{itemize}

The first was proved by Bestvina~\cite[Theorem 4.5]{Bes99}
and the second by
Charney, Meier and Whittlesey~\cite[Corollary 6.9]{CMW04}
following the arguments of Bestvina.
For the full statement of their results, see \S\ref{ssec:CentQuot}.

We show that Bestvina's result holds for all Garside groups.

\medskip\noindent
\textbf{Theorem~\ref{thm:cyclic}.}\ \ \em
Let $G$ be a Garside group with Garside element $\Delta$.
Then every finite subgroup of\/ $G_\Delta$ is cyclic.
\upshape
\medskip

Our proof uses the result of Charney, Meier and Whittlesey.
Actually we prove that every finite abelian subgroup of $G_\Delta$ is cyclic.
Because every finite subgroup of $G_\Delta$ is abelian,
this implies the above theorem.

\subsection{Primitive periodic elements}

Let us recall the braid group $B_n$, which is the same as the Artin group
$A(\arA_{n-1})$.
It has the group presentation~\cite{Art25}:
$$
B_n  =  \left\langle \sigma_1,\ldots,\sigma_{n-1} \left|
\begin{array}{ll}
\sigma_i \sigma_j = \sigma_j \sigma_i & \mbox{if } |i-j| > 1, \\
\sigma_i \sigma_j \sigma_i = \sigma_j \sigma_i \sigma_j & \mbox{if } |i-j| = 1.
\end{array}
\right.\right\rangle.
$$
The group $B_n$ admits two well-known Garside structures:
the classical Garside structure~\cite{Gar69, EC+92, EM94}
and the dual Garside structure~\cite{BKL98}.
Let
$\Delta = \sigma_1(\sigma_2\sigma_1)\cdots(\sigma_{n-1}\cdots\sigma_2\sigma_1)$ and
$\delta = \sigma_{n-1}\sigma_{n-2}\cdots\sigma_1$.
Then $\Delta$ and $\delta$ are the Garside elements in the classical and
dual Garside structures, respectively.

\smallskip
The most fundamental question on periodic elements in Garside groups
would be a characterization of them.
For the braid group $B_n$, it is a classical theorem of
Brouwer, Ker\'ekj\'art\'o and Eilenberg~\cite{Bro19, Ker19, Eil34,BDM02}
that an $n$-braid is periodic if and only if
it is conjugate to a power of either $\delta$ or $\varepsilon$,
where $\varepsilon=\delta\sigma_1$.
The same kind of statement holds
for the Artin groups $A(\arB_n)$, $A(\arD_n)$, $A(\arI_2(e))$ and
the braid group $B(e,e,n)$ of the complex reflection group of type $(e,e,n)$
(see Theorem~\ref{thm:per-elt}).
It is worth mentioning that Bessis~\cite{Bes06a, Bes06b}
explored many important properties
of periodic elements
in the context of braid groups of complex reflection groups.

For arbitrary Garside groups, we cannot expect such a nice characterization
of periodic elements.
So we establish a weaker theorem.
Let $G$ be a Garside group with Garside element $\Delta$.
Let us say a nonidentity element $g\in G$ to be \emph{primitive}
if it is not a nontrivial power of another element,
that is, $g=h^k$ for $h\in G$ and $k\in\Z$ implies $k=\pm1$.
Then the Brower-Ker\'ekj\'art\'o-Eilenberg theorem can be restated as:
the braids $\delta$ and $\varepsilon$ are the only primitive
periodic elements, up to conjugacy and taking inverse.
Since $\delta^n=\varepsilon^{n-1}=\Delta^2$ and $\Delta^2$ is central,
every primitive periodic braid is a root of $\Delta^2$.
We generalize this property to arbitrary Garside groups.

\medskip\noindent
\textbf{Theorem~\ref{thm:primitive}.}\ \ \em
Let $G$ be a Garside group with Garside element $\Delta$,
and $\Delta^m$ the minimal positive central power of $\Delta$.
Then every primitive periodic element in $G$ is a $k$-th root of $\Delta^m$
for some $k$ with $1\le | k | \le m\Vert\Delta\Vert$.
\upshape
\medskip

Using the above theorem, we show in Proposition~\ref{prop:root} that
there is a finite-time algorithm that,
given a Garside group, computes
all primitive periodic elements up to conjugacy.
Therefore, in theory, though probably difficult in practice,
we can establish a Brower-Ker\'ekj\'art\'o-Eilenberg type theorem
for any fixed Garside group.

\subsection{Conjugacy classes of periodic elements}\label{ssec:ConjClass}

The conjugacy problem in a group has two versions:
the conjugacy decision problem (CDP) is to decide
whether given two elements are conjugate or not;
the conjugacy search problem (CSP) is to find a conjugating element
for a given pair of conjugate elements.
In the late sixties Garside~\cite{Gar69} first solved
the conjugacy problem in braid groups.
Then there have been considerable efforts to improve
his solution~\cite{BS72,Del72,EC+92,EM94}.

The CDP for periodic braids is easy:
an $n$-braid $\alpha$ is periodic if and only if
either $\alpha^n$ or $\alpha^{n-1}$ belongs to the cyclic group $\langle\Delta^2\rangle$;
two periodic braids are conjugate if and only if
they have the same exponent sum.
The situation is similar for the groups
$A(\arB_n)$, $A(\arD_n)$, $A(\arI_2(e))$ and $B(e,e,n)$.
In these groups, there is an easy periodicity test for elements
(see \S\ref{ssec:periodicity}),
and two periodic elements are conjugate if and only if
they have the same exponent sum (see Proposition~\ref{prop:ExpSumConj}).

The CSP for periodic braids is not as easy as the CDP.
The standard solution is not efficient enough.
Recently Birman, Gebhardt and Gonz\'alez-Meneses~\cite{BGG07}
constructed an efficient solution,
by using several known isomorphisms between Garside structures
on the braid groups and other Garside groups.
However, unlike the case of CDP, their solution does not naively
extend to other Garside groups such as $A(\arD_n)$ and $B(e,e,n)$,
because their isomorphisms are peculiar to braids.

For the CDP and CSP in arbitrary Garside groups, as far as the authors know,
there is no solution specialized to periodic elements.
On this account, we study properties of conjugacy classes
of periodic elements
in arbitrary Garside groups.

\smallskip

Let $G$ be a Garside group with Garside element $\Delta$,
and $\Delta^m$ the minimal positive central power of $\Delta$.
If $g\in G$ is $p/q$-periodic, $g^q$ is conjugate to $\Delta^p$ and
$q$ is the smallest among positive integers with such property.
We define a $p/q$-periodic element to be \emph{precentral} if $p\equiv 0\bmod m$, and
\emph{slim} if $p\equiv 1\bmod q$.

Let $g=\Delta^u a_1a_2\cdots a_\ell \in G$ be in normal form,
and let $b$ be a prefix of $a_1$.
The conjugation
$$
\tau^{-u}(b)^{-1} g \tau^{-u}(b) = \Delta^u a_1'a_2\cdots a_\ell\tau^{-u}(b)
$$
is called a \emph{partial cycling} of $g$ by $b$,
where $\tau(x)=\Delta^{-1}x\Delta$ and $a_1'=b^{-1}a_1$.

We establish the following theorem, where
$[g]^{\inf}$, $[g]^S$, $[g]^U$ and $[g]^{St}$ denote the summit set,
super summit set, ultra summit set and
stable super summit set of $g$, respectively.

\medskip\noindent
\textbf{Theorem~\ref{thm:pa-cy}.}\ \ \em
Let $g$ be a slim, precentral periodic element
of a Garside group $G$.
Then
$$
[g]^{\inf}=[g]^S=[g]^U=[g]^{St}.
$$
In particular, $[g]^{S}$ is closed under any partial cycling.
\upshape
\medskip

The above theorem will be useful in solving
the CSP for periodic elements, at least in the groups
$A(\arA_n)$, $A(\arD_n)$ and $B(e,e,n)$.
The complexity of the standard conjugacy algorithm is proportional to
the size of the super summit set,
which is exponential with respect to the braid index
in the case of braid groups
as observed by Birman, Gebhardt and Gonz\'alez-Meneses~\cite{BGG07}.
Therefore, in order to make an efficient algorithm,
we need further information on super summit sets, like Theorem~\ref{thm:pa-cy}.
We hope that our work will be also useful in studying
conjugacy classes of periodic elements in other Garside groups.

In Theorem~\ref{thm:pa-cy}, a periodic element is required to be
both slim and precentral.
These requirements are necessary (see Example~\ref{eg:conditions}).
In solving the conjugacy problem for periodic elements in Garside groups,
we may assume without loss of generality that given periodic elements are slim
(see Lemma~\ref{lem:BCMW}).
As for the precentrality condition,
we can make every periodic element precentral by modifying the Garside structure:
for a Garside group $G$ with Garside element $\Delta$,
if we change the Garside structure on $G$ by declaring
the central power $\Delta^m$ as a new Garside element, then every periodic element
becomes precentral in this new Garside structure.
When working with the groups $A(\arA_n)$, $A(\arB_n)$, $A(\arD_n)$, $A(\arI_2(e))$ and $B(e,e,n)$
endowed with the dual Garside structure,
we may further assume the precentrality by the following theorem.

\medskip\noindent
\textbf{Theorem~\ref{thm:dual}.}\ \ \em
In the dual Garside structure on each of the groups
$A(\arA_n)$, $A(\arB_n)$, $A(\arD_n)$, $A(\arI_2(e))$ and $B(e,e,n)$,
every periodic element is either precentral
or conjugate to a power of the Garside element.
\upshape

\subsection{Organization}

Section 2 provides a brief introduction to Garside groups.
Section 3 studies periodic elements in Garside groups.
Section 4 studies periodic elements in the groups
$A(\arA_n)$, $A(\arB_n)$, $A(\arD_n)$, $A(\arI_2(e))$ and $B(e,e,n)$.
Section 5 discusses some algorithmic problems concerning
periodic elements, such as periodicity decision problem,
tabulation of all primitive periodic elements,
and so on.

\section{Review of Garside groups}

\subsection{Garside groups}

Garside groups were defined by Dehornoy and Paris~\cite{DP99}
as groups satisfying certain conditions under which the strategy and results of
Garside~\cite{Gar69}, Deligne~\cite{Del72}, Brieskorn and Saito~\cite{BS72} still hold.
This section briefly reviews Garside groups.
For a detailed description, see~\cite{DP99,Deh02}.

\begin{definition}
For a monoid $M$, let $1$ denote the identity element.
An element $a\in M\setminus \{ 1\}$ is called an \emph{atom} if
$a=bc$ for $b,c\in M$ implies either $b=1$ or $c=1$.
For $a\in M$, let $\Vert a\Vert$ be the supremum
of the lengths of all expressions of
$a$ in terms of atoms. The monoid $M$ is said to be \emph{atomic}
if it is generated by its atoms and $\Vert a\Vert<\infty$
for any element $a$ of $M$.
In an atomic monoid $M$, there are partial orders $\le_L$ and $\le_R$:
$a\le_L b$ if $ac=b$ for some $c\in M$;
$a\le_R b$ if $ca=b$ for some $c\in M$.
\end{definition}

\begin{definition}
An atomic monoid $M$ is called a \emph{Garside monoid} if it satisfies the following.
\begin{enumerate}
\item
$M$ is left and right cancellative.

\item
$(M,\le_L)$ and $(M,\le_R)$ are lattices.
That is, for every $a, b\in M$ there are a unique least common multiple $a\vee_L b$ (resp.{} $a\vee_R b$)
and a unique greatest common divisor $a\wedgeL b$ (resp.{} $a\wedgeR b$) with respect to
$\le_L$ (resp.{} $\le_R$).

\item
$M$ contains an element $\Delta$, called a
\emph{Garside element}, satisfying the following:
\begin{enumerate}
\item
For each $a\in M$, $a\le_L\Delta$ if and only if $a\le_R\Delta$.

\item
The set
$[1, \Delta] = \{ a\in M \mid 1\le_L a\le_L\Delta\}$ is finite and generates $M$.
Elements of this set are called {\em simple elements}.
$(1,\Delta)$ denotes the set $[1,\Delta]\setminus\{1,\Delta\}$.
Similarly for $[1,\Delta)$ and $(1,\Delta]$.
\end{enumerate}
\end{enumerate}
\end{definition}

\begin{definition}
Let $M$ be a Garside monoid with Garside element $\Delta$.
The group $G$ of fractions of $M$ is called a \emph{Garside group}.
We identify the elements of $M$ and their images in $G$,
and call them \emph{positive elements} of $G$.
The Garside monoid $M$ is often denoted by $G^+$.
We will call the pair $(G^+, \Delta)$ a {\em Garside structure} on $G$.
\end{definition}

Note that a Garside group can have more than one Garside structure.
From now on, $G$ denotes a Garside group
with a fixed Garside structure $(G^+, \Delta)$, if not specified otherwise.

\begin{definition}
The partial orders $\le_L$ and $\le_R$, and thus the lattice structures
in the Garside monoid $G^+$ can be extended to the Garside group $G$.
For $a, b\in G$, $a\le_L b$ (resp. $a\le_R b$) means $a^{-1}b\in G^+$
(resp. $ba^{-1}\in G^+$).
\end{definition}

\begin{definition}
Let $\tau : G\to G$ be the inner automorphism of $G$
defined by $\tau(g)=\Delta^{-1}g\Delta$ for $g\in G$.
\end{definition}

The automorphism $\tau$ preserves
the set $[1,\Delta]$ which is a finite set of generators of $G$.
Therefore some power of $\tau$ is the identity,
equivalently, some power of $\Delta$ is central.

\begin{definition}
For $a,b\in G^+$, $a$ is called a \emph{prefix} (resp.\ \emph{suffix}) of $b$
if $a\le_L b$ (resp.\ $a\le_R b$).
\end{definition}

\begin{definition}
For every $g\in G$, there exists a unique decomposition
$$g=\Delta^u a_1\cdots a_\ell$$
such that $u\in\Z$, $\ell\in\Z_{\ge 0}$, $a_1,\ldots,a_\ell\in (1,\Delta)$ and $(a_i a_{i+1}\cdots a_\ell)\wedgeL \Delta=a_i$
for $i=1,\ldots,\ell$.
This decomposition is called the \emph{(left) normal form} of $g$.
In this case, $\inf(g)=u$, $\sup(g)=u+\ell$ and $\len(g)=\ell$ are
called the \emph{infimum}, \emph{supremum} and \emph{canonical length} of $g$,
respectively.
\end{definition}

\begin{definition}
Let $g=\Delta^u a_1\cdots a_\ell\in G$ be in normal form with $\ell\ge 1$.
The {\em cycling} $\c(g)$ and {\em decycling} $\d(g)$ of $g$ are conjugations of $g$
defined as
\begin{align*}
\c(g)&= \Delta^u a_2\cdots a_\ell\tau^{-u}(a_1) = (\tau^{-u}(a_1))^{-1}g \tau^{-u}(a_1),\\
\d(g)&= \Delta^u \tau^{u}(a_\ell)a_1\cdots a_{\ell-1} = a_{\ell} g a_{\ell}^{-1}.
\end{align*}
We define $\c(\Delta^u)=\d(\Delta^u)=\Delta^u$ for $u\in\Z$.
\end{definition}

\begin{definition}
For $g\in G$, we denote its conjugacy class $\{ h^{-1}gh : h\in G\}$ by $[g]$.
The conjugacy invariants $\infs(g)$, $\sups(g)$ and $\lens(g)$ are defined as follows:
$\infs(g)=\max\{\inf(h):h\in [g]\}$;
$\sups(g)=\min\{\sup(h):h\in [g]\}$;
$\lens(g)=\sups(g)-\infs(g)$.
They are called the
\emph{summit infimum},
\emph{summit supremum} and
\emph{summit canonical length}
of $g$, respectively.
The \emph{summit set} $[g]^{\inf}$,
\emph{super summit set} $[g]^S$,
\emph{ultra summit set} $[g]^U$ and
\emph{stable super summit set} $[g]^{St}$
are defined as follows:
\begin{align*}
[g]^{\inf} &=\{h\in[g]: \inf(h)=\infs(g)\};\\{}
[g]^S&=\{h\in [g]:\inf(h)=\infs(g) \ \mbox{ and } \sup(h)=\sups(g)\};\\{}
[g]^{U}&=\{h\in [g]^S: \c^k(h)=h \ \mbox{ for some positive integer $k$}\};\\{}
[g]^{St}&=\{h\in [g]^S:h^k\in[g^k]^S \ \mbox{ for all positive integers $k$}\}.
\end{align*}
\end{definition}

For every $g\in G$, the sets $[g]^{\inf}$, $[g]^S$, $[g]^U$ and
$[g]^{St}$ are all finite, nonempty and computable in a finite
number of steps~\cite{EM94,Geb05,LL08a}.

\subsection{Translation discreteness}
For every $g$ in a Garside group $G$, the following limits are well-defined:
$$
\INF(g)=\lim_{n\to\infty}\frac{\inf(g^n)}n;\quad
\SUP(g)=\lim_{n\to\infty}\frac{\sup(g^n)}n;\quad
\LEN(g)=\lim_{n\to\infty}\frac{\len(g^n)}n.
$$
Observe that $\SUP(g)=\INF(g)+\LEN(g)$.
These limits were defined in order to investigate discreteness of
translation numbers in Garside groups~\cite{LL07,LL08b}.
They can be computed explicitly,
for example by using the formulas~\cite[Theorem 5.1]{LL08b}
$$\begin{array}{rcl}
\INF(g)&=&\max\{\infs(g^k)/k : k=1,\ldots, \Vert\Delta\Vert\},\\
\SUP(g)&=&\min\{\sups(g^k)/k : k=1,\ldots, \Vert\Delta\Vert\}.
\end{array}
$$

For a group $G$ with a set of generators $X$,
the translation number of $g\in G$ with respect to $X$
is defined by
$t_X(g)=\lim_{n\to\infty} \frac{|g^n|_X}n$,
where $|g^n|_X$ denotes
the minimal word-length of $g^n$ with respect to $X$.
In a Garside group $G$ with Garside element $\Delta$
the above limits are related to the translation number by
$$
t_X(g)=\max\{~ |\INF(g)|,~ |\SUP(g)|,~ |\LEN(g)|~\}
$$
provided the set $[1,\Delta]$ is taken as $X$.

The following proposition collects some important properties
of the above limits.

\begin{proposition}[{\cite{LL07,LL08b}}]\label{prop:BasicTinf}
For $g$ and $h$ in a Garside group $G$ with Garside element $\Delta$,
\begin{enumerate}
\item
$\INF(h^{-1}gh)=\INF(g)$ and $\SUP(h^{-1}gh)=\SUP(g)$;

\item
$\INF(g^n)= n\cdot\INF(g)$ and $\SUP(g^n)= n\cdot\SUP(g)$
for all integers $n\ge 1$;

\item
$\infs(g)=\lfloor \INF(g)\rfloor$ and $\sups(g)=\lceil \SUP(g)\rceil$;

\item
$\INF(g)$ and $\SUP(g)$ are rational of the form $p/q$,
where $p$ and $q$ are coprime integers and
$1\le q\le\Vert\Delta\Vert$.
\end{enumerate}
\end{proposition}

\section{Periodic elements in Garside groups}

If not specified otherwise,
$G$ is a Garside group with Garside element $\Delta$,
and $\Delta^m$ is the minimal positive central power of $\Delta$.
This section studies periodic elements in Garside groups.

\subsection{Periodicity and Garside structure}\label{ssec:Garsid structure}

For a group $G$, let $Z(G)$ denote the center of $G$.

\begin{theorem}\label{thm:Gar_struc}
Let $G$ be a Garside group.
Then $Z(G)$ is cyclic if and only if
any pair of Garside elements of $G$ are commensurable.
\end{theorem}

It is easy to see that
the periodicity of an element in $G$ does not depend on
the choice of a particular Garside element
if and only if any pair of Garside elements of $G$ are commensurable.
Hence we have the following.

\begin{corollary}\label{cor:Gar_struc}
The periodicity of an element of\/ $G$ does not depend on the choice of a particular
Garside element if and only if $Z(G)$ is cyclic.
\end{corollary}

\smallskip
We prove Theorem~\ref{thm:Gar_struc} by using the following lemma.

\begin{lemma}\label{lem:Gar_struc}
Let $(G^+, \Delta)$ be a Garside structure on a group $G$.
For $g\in G$, let
$L(g)=\{ a\in G^{+} : a\le_L g \}$ and $R(g)=\{ a\in G^{+} : a\le_R g \}$.
\begin{itemize}
\item[(i)]
Let $c$ be a positive element in $Z(G)$.
Then $L(c)=R(c)$.

\item[(ii)]
Let $c$ be a positive element in $Z(G)$ with $\Delta\le_L c$.
Then $c$ is a Garside element, that is,
$L(c)=R(c)$ and $L(c)$ generates the Garside monoid $G^+$.
\end{itemize}
\end{lemma}

\begin{proof}
(i)\ \
Let $a\in L(c)$, then $c=ab$ for some $b\in G^+$.
Because $c$ is central, $ab=c=bcb^{-1}=b(ab)b^{-1}=ba$.
Therefore $c=ba$, hence $a\in R(c)$.
This means that $L(c)\subset R(c)$.
Similarly, $R(c)\subset L(c)$.

\smallskip(ii)\ \
By (i), $L(c)=R(c)$.
As $\Delta\le_L c$, we have $L(\Delta)\subset L(c)$.
Since $L(\Delta)$ generates $G^+$, so does $L(c)$.
\end{proof}

\begin{proof}[Proof of Theorem~\ref{thm:Gar_struc}]
Suppose that $Z(G)$ is cyclic.
Let $(G_1^+, \Delta_1)$ and $(G_2^+, \Delta_2)$ be Garside structures on $G$.
Then there exist positive integers $m_1$ and $m_2$ such that
$\Delta_1^{m_1}$ and $\Delta_2^{m_2}$ are central in $G$.
Because $Z(G)$ is cyclic,
$\Delta_1^{m_1}$ and $\Delta_2^{m_2}$ are commensurable,
hence $\Delta_1$ and $\Delta_2$ are commensurable.

\smallskip
Conversely, suppose any pair of Garside elements are commensurable.
Fix a Garside structure $(G^+, \Delta)$ on $G$.
Let $m$ be the smallest positive integer such that $\Delta^m$ is central.

We claim that any nonidentity element of $Z(G)$ is commensurable with $\Delta$.
Let $g$ be a nonidentity central element.
Take an integer $k$ such that
$$
k\equiv 0\bmod m\quad\mbox{and}\quad k\ge -\inf(g)+1.
$$
Let $c=\Delta^k g$. Then $c$ is a central element with $\Delta\le_L c$,
hence $c$ is a Garside element by Lemma~\ref{lem:Gar_struc}.
By the hypothesis, $c$ is commensurable with $\Delta$.
As $c$ is central, there exist nonzero integers $p$ and $q$ such that
$\Delta^p=c^q=(\Delta^kg)^{q} = \Delta^{kq}g^q$.
Since $g^q=\Delta^{p-kq}$, $g$ is commensurable with $\Delta$.

It is known that Garside groups are torsion-free
by Dehornoy~\cite{Deh98},
and that every abelian subgroup of a Garside group is
finitely generated by Charney, Meier and Whittlesey~\cite{CMW04}.
Thus $Z(G)$ is torsion-free and finitely generated.
Moreover, by the above claim, any two nonidentity elements of $Z(G)$ are commensurable
because each of them is commensurable with $\Delta$.
These imply that $Z(G)$ is cyclic.
\end{proof}

The following example shows that periodicity of an element
depends on the choice of a particular Garside element.
\begin{example}\label{ex:Per_GarStru}
The group $\Z\times\Z$ is a Garside group with Garside monoid $\Z_{\ge 0}\times\Z_{\ge 0}$.
For any $a,b\in\Z_{\ge 1}$, the element $(a,b)$ is a Garside element.
In particular, $\Delta_1=(2,2)$ and $\Delta_2=(2,3)$ are Garside elements.
Let $g=(1,1)$. Then $g$ is periodic with respect to $\Delta_1$
but not to $\Delta_2$.
\end{example}

\subsection{Roots of periodic elements}\label{ssec:roots}

Periodic elements can be defined as follows.

\begin{definition}\label{defn:pq-per}
An element $g\in G$ is said to be \emph{$p/q$-periodic}
for $p\in\Z$ and $q\in\Z_{\ge 1}$ if $g^q$ is conjugate to $\Delta^p$ and
$q$ is the smallest positive integer such that $g^q$
is conjugate to a power of $\Delta$.
\end{definition}

The following lemma shows some properties of periodic elements regarding translation numbers.
For integers $p$ and $q$ with at least one of them different from zero,
$p\wedge q$ denotes their greatest common divisor.

\begin{lemma}\label{lem:per-elt}
Let $g\in G$ be a periodic element. Then the following hold.
\begin{enumerate}
\item
$\LEN(g)=0$, that is, $\INF(g)=\SUP(g)$.

\item
$\INF(g^k) = k\cdot\INF(g)$ for all $k\in\Z$.

\item
For any $k\in\Z$, $\lens(g^k)$ is either 0 or 1.

\item
Let $\INF(g)=p/q$ for $p\in\Z$ and $q\in\Z_{\ge 1}$ with $p\wedge q=1$.
Then the following are equivalent for $k\in\Z$:
\begin{enumerate}
\item
$g^k$ is conjugate to a power of $\Delta$;

\item
$\lens(g^k)=0$;

\item
$\INF(g^k)$ is an integer;

\item
$k$ is a multiple of $q$.
\end{enumerate}
In particular, $g^q$ is conjugate to $\Delta^p$.
\end{enumerate}
\end{lemma}

\begin{proof}
(i)\ \
Since $g^k=\Delta^\ell$ for some $k\in\Z_{\ge 1}$ and $\ell\in\Z$,
$$
\LEN(g)=\frac1k\cdot\LEN(g^k)=\frac1k\cdot\LEN(\Delta^\ell)
=\frac1k\cdot 0=0.
$$

\smallskip(ii)\ \
We know that $\INF(g^k) = k\cdot\INF(g)$ holds for all $k\ge 0$.
Let $k<0$, then $k=-\ell$ for some $\ell\ge 1$.
Since $\INF(g)=\SUP(g)$ by (i)
and $\INF(h^{-1})=-\SUP(h)$ for all $h\in G$,
we have
$$
\INF(g^k) =\INF((g^{-1})^\ell)
=\ell\cdot\INF(g^{-1})
=\ell\cdot(-\SUP(g))=k\cdot\INF(g).
$$

\smallskip(iii) \ \
Let $\INF(g)=p/q$ for $p\in\Z$ and $q\in\Z_{\ge 1}$ with $p\wedge q=1$.
Choose any $k\in\Z$.
Because $\SUP(g)=\INF(g)=p/q$ by (i),
\begin{equation}
\label{eq:lens}
\lens(g^k)
=\sups(g^k)-\infs(g^k)
=\lceil \SUP(g^k)\rceil -\lfloor\INF(g^k)\rfloor
=\lceil kp/q\rceil -\lfloor kp/q\rfloor
\end{equation}
by (ii). Therefore $\lens(g^k)$ is either 0 or 1.

\smallskip (iv) \ \
It is obvious that $\lens(g^k)=0$ if and only if $g^k$ is conjugate
to a power of $\Delta$.
By Eq.~(\ref{eq:lens}),
$\lens(g^k)=0$
if and only if $\INF(g^k)=kp/q$ is an integer.
Because $p$ and $q$ are relatively prime,
$kp/q$ is an integer if and only if $k$ is a multiple of $q$.
Therefore the four conditions---(a), (b), (c) and (d)---are equivalent.

By (a) and (d), $g^q$ is conjugate to $\Delta^\ell$ for some
integer $\ell$, hence $\INF(g)=\ell /q$.
Because $\INF(g)=p/q$ by the hypothesis, we have $\ell=p$.
\end{proof}

\begin{remark}
If $g$ is an element of $G$ with $\LEN(g)=0$, then $\INF(g)=\SUP(g)=p/q$
for some $p\in\Z$ and $q\in\Z_{\ge 1}$.
Hence $\infs(g^q)=\lfloor \INF(g^q)\rfloor=p=\lceil \SUP(g^q)\rceil=\sups(g^q)$.
This implies that $g^q$ is conjugate to $\Delta^p$,
hence $g$ is periodic.
Combining with Lemma~\ref{lem:per-elt}(i), we can see that $g$ is periodic if and only if $\LEN(g)=0$.
\end{remark}

The following corollary is used in later sections.
\begin{corollary}\label{cor:basic}
Let $g\in G$ be periodic, and let $k$ be an integer.
\begin{enumerate}
\item
$\INF(g)=k$\/ if and only if $g$ is conjugate to $\Delta^k$.
\item
$\INF(g)=mk$\/ if and only if $g=\Delta^{mk}$.
\end{enumerate}
\end{corollary}

\begin{proof}
If $\INF(g)=k$, then $g$ is conjugate to $\Delta^k$ by Lemma~\ref{lem:per-elt}(iv).
The converse direction is obvious.
This proves (i), and
(ii) follows immediately from (i) as $\Delta^{mk}$ is central.
\end{proof}

\begin{theorem}\label{thm:unique}
Let $G$ be a Garside group with Garside element $\Delta$, and let $g\in G$
and $a, b, k \in\Z_{\neq 0}$.
\begin{itemize}
\item[(i)]
If\/ $g^{kb}$ is conjugate to $\Delta^{ka}$,
then $g^b$ is conjugate to $\Delta^a$.

\item[(ii)]
If\/ each of $g^{a}$ and $g^{b}$ is conjugate to a power of $\Delta$,
then so is $g^{a\wedge b}$.
\end{itemize}
\end{theorem}

\begin{proof}
The hypothesis in either case of (i) or (ii) implies that $g$ is periodic.
Let $\INF(g)=p/q$ for $p\in\Z$ and $q\in\Z_{\ge 1}$ with $p\wedge q=1$.
Then $g^{q}$ is conjugate to $\Delta^{p}$ by Lemma~\ref{lem:per-elt}.

\smallskip

(i)\ \
Since $g^{kb}$ is conjugate to $\Delta^{ka}$, one has
$\INF(g)=a/b=p/q$, hence there is $d\in\Z_{\neq 0}$ such that $a=dp$ and $b=dq$.
Therefore $g^b$ is conjugate to $\Delta^a$.

\smallskip

(ii)\ \
By Lemma~\ref{lem:per-elt}, both $a$ and $b$ are multiples of $q$,
hence $a\wedge b$ is a multiple of $q$.
Therefore $g^{a\wedge b}$ is conjugate to some power of\/ $\Delta$.
\end{proof}

\begin{corollary}\label{cor:Periodic_Tinf}
Let $g\in G$ be a periodic element.
Then, $g$ is $p/q$-periodic if and only if
$\INF(g)=p/q$ for $p\in\Z$ and $q\in\Z_{\ge 1}$ with $p\wedge q=1$.
\end{corollary}

\begin{proof}
Suppose that $\INF(g)=p/q$ for $p\in\Z$ and $q\in\Z_{\ge 1}$ with $p\wedge q=1$.
Then, by Lemma~\ref{lem:per-elt}(iv), $g^q$ is conjugate to $\Delta^p$.
Assume that $g^b$ is conjugate to $\Delta^a$ for $a, b\in\Z$ with $1\le b<q$.
Then $\INF(g)=a/b = p/q$, which implies $q\le b$.
It is a contradiction.
Thus $g$ is $p/q$-periodic.

The converse direction is obvious by Theorem~\ref{thm:unique}.
\end{proof}

Theorem~\ref{thm:unique}(i) implies that
if $g^k=\Delta^{k\ell}$ for a nonzero integer $k$
then $g$ is conjugate to $\Delta^\ell$.
The following example illustrates that
the general statement for the uniqueness of roots up to conjugacy
(i.e.\ if $g^k = h^k$ for $k\neq 0$ then $g$ is conjugate to $h$)
does not hold in Garside groups, even for periodic elements.

\begin{example}\label{ex:nonunique}
Let $G$ be the group defined by
$$
G=\langle x,y\mid x^a=y^a\rangle,\qquad a\ge 2.
$$
It is a Garside group with Garside element
$\Delta=x^a=y^a$~\cite[Example 4]{DP99}.
(Note that $x$ and $y$ are periodic elements and that $\Delta$ is central.)
We claim that
\begin{itemize}
\item[(i)] $x$ and $y$ are not conjugate;
\item[(ii)] there is no Garside structure on $G$
in which $x$ is a Garside element.
\end{itemize}

Because $G/\langle \Delta\rangle =\langle x,y\mid x^a=y^a=1\rangle
=\langle x\mid x^a=1\rangle*\langle y\mid y^a=1\rangle$,
the images of $x$ and $y$ in $G/\langle \Delta\rangle$ are not conjugate.
Therefore $x$ and $y$ are not conjugate.

Assume that there exists a Garside structure on $G$ with Garside element $x$.
Because $x^a=y^a$ and $x$ is a Garside element,
$y$ is conjugate to $x$ by Theorem~\ref{thm:unique}(i).
It is a contradiction to (i).
\end{example}

The above example shows that
\begin{quote}
there is a Garside group with a periodic element $g$
such that there is no Garside structure in which
$g$ is a Garside element.
\end{quote}
Therefore it gives a negative answer to the question of Bessis
stated in \S\ref{ssec:Intr_Root}.

\subsection{Primitive periodic elements}\label{ssec:Primitive}

The famous theorem of Brower, Ker\'ekj\'art\'o and Eilenberg~\cite{Bro19,Ker19,Eil34} says that
in the braid group $B_n$, there are two periodic elements $\delta$ and $\varepsilon$
such that every other periodic element is conjugate to
a power of either $\delta$ or $\varepsilon$.
Motivated by this, we define the following notion.

\begin{definition}
A nonidentity element $g\in G$ is said to be \emph{primitive}
if it is not a nontrivial power of another element,
that is, $g=h^k$ for $h\in G$ and $k\in\Z$ implies $k=\pm1$.
\end{definition}

Using the above terminology
the Brower-Ker\'ekj\'art\'o-Eilenberg theorem
can be restated as:
the braids $\delta$ and $\varepsilon$ are the only primitive
periodic elements in the braid group $B_n$, up to inverse and conjugacy.
Notice that both $\delta$ and $\varepsilon$ are roots of $\Delta^2$
as $\delta^n=\Delta^2=\varepsilon^{n-1}$.
Therefore every primitive periodic braid is a root of $\Delta^2$.
We generalize this property to Garside groups in Theorem~\ref{thm:primitive}.
The following lemma is a key to doing this,
and will be used later as well.

For a $p/q$-periodic element $g\in G$, $p=0$ if and only if $g$ is the identity,
because Garside groups are torsion-free~\cite{Deh98}.

\begin{lemma}\label{lem:PerCyclic}
Let $g\in G$ be $p/q$-periodic with $p\neq 0$,
and let $H$ be the subgroup of\/ $G$ generated by $g$ and $\Delta^m$.
Let $h=g^{r}\Delta^{ms}$, where
$r$ and $s$ are integers  with $pr + qms= p\wedge m$.
Then $H$ is a cyclic group generated by $h$.
More precisely, $g = h^{\frac{p}{p\wedge m}}$ and $\Delta^m = h^{\frac{qm}{p\wedge m}}$.
\end{lemma}

\begin{proof}
Since $p$ is coprime to $q$, one has $p\wedge qm = p\wedge m$,
hence there are integers $r$ and $s$ with
$pr + qms= p\wedge m$.
As $g^q$ is conjugate to $\Delta^p$ and $\Delta^m$ is central,
one has $g^{q\frac{m}{p\wedge m}}=\Delta^{p\frac{m}{p\wedge m}}$.
Using this identity, we have
\begin{align*}
h^{\frac{p}{p\wedge m}} & = (g^{r}\Delta^{ms})^{\frac{p}{p\wedge m}}
  = g^{\frac{pr}{p\wedge m}} \Delta^{p\frac{m}{p\wedge m}s}
  = g^{\frac{pr}{p\wedge m}} g^{q\frac{m}{p\wedge m}s}
  = g^{\frac{pr+qms}{p\wedge m}} = g,\\
h^{\frac{qm}{p\wedge m}} & = (g^{r}\Delta^{ms})^{\frac{qm}{p\wedge m}}
  = g^{q\frac{m}{p\wedge m}r} \Delta^{\frac{qms}{p\wedge m}m}
  = \Delta^{p\frac{m}{p\wedge m}r} \Delta^{\frac{qms}{p\wedge m}m}
  = \Delta^{m\frac{pr+qms}{p\wedge m}} = \Delta^m.
\end{align*}
Therefore $H$ is a cyclic group generated by $h$.
\end{proof}

\begin{theorem}\label{thm:primitive}
Every primitive periodic element in $G$ is a $k$-th root of $\Delta^m$
for some $k$ with $1\le | k | \le m\Vert\Delta\Vert$.
\end{theorem}

\begin{proof}
Let $g$ be a primitive $p/q$-periodic element of $G$.
Since $g=h^{\frac{p}{p\wedge m}}$ for some $h\in G$
by Lemma~\ref{lem:PerCyclic},
we have $\frac{p}{p\wedge m}=\pm 1$,
hence $m$ is a multiple of $p$.
So $\frac{qm}{p}$ is an integer.
As $\INF(g^{\frac{qm}{p}})=\frac{qm}{p}\cdot\frac pq=m$,
we have $g^{\frac{qm}{p}}=\Delta^m$ by Corollary~\ref{cor:basic}(ii).
Therefore $g$ is a $\frac{qm}{p}$-th root of $\Delta^m$.
By Proposition~\ref{prop:BasicTinf}, we know $1\le q\le \Vert\Delta\Vert$.
Therefore $1\le \left|\frac{qm}{p}\right| = \left|\frac{m}{p}\right| \cdot q \le m\Vert\Delta\Vert$.
\end{proof}

\subsection{Quotient group $G_\Delta$}
\label{ssec:CentQuot}

Let $G_\Delta$ be the quotient $G/\langle\Delta^m\rangle$,
where $\myangle{\Delta^m}$ is the cyclic group generated by $\Delta^m$.
For an element $g\in G$, let $\bar g$ denote the image of $g$
under the natural projection from $G$ to $G_\Delta$.
Hence,  $g\in G$ is periodic if and only if
$\bar g$ is of finite order in $G_\Delta$.
For periodic elements in $G$, it is sometimes more convenient to view them in $G_\Delta$.

The following theorem was proved
by Bestvina~\cite[Theorem 4.5]{Bes99} for
Artin groups of finite type, and then
proved by Charney, Meier and Whittlesey~\cite[Corollary 6.9]{CMW04}
for Garside groups.

\begin{theorem}[\cite{Bes99,CMW04}]
\label{thm:CMW}
The finite subgroups of\/ $G_\Delta$ are, up to conjugacy,
one of the following two types:
\begin{itemize}
\item[(i)]
the cyclic group generated by the image of $\Delta^u a$ in $G_\Delta$
for some $u\in\Z$ and some simple element $a\ne\Delta$
such that if $a\neq 1$, then for some integer $2\le q\le\Vert\Delta\Vert$
$$
\tau^{(q-1)u}(a)\,\tau^{(q-2)u}(a)\cdots \tau^{u}(a)\,a=\Delta ;
$$

\item[(ii)]
the direct product of a cyclic group of type (i) and
$\langle\bar\Delta^k\rangle$ where $\Delta^k$ commutes with $a$.
\end{itemize}
\end{theorem}

In the case of Artin groups of finite type,
Bestvina showed that finite subgroups of $G_\Delta$ are all cyclic groups
(hence they are of type (i) in the above theorem).
Using the following lemma, we show in Theorem~\ref{thm:cyclic}
that the same is true for Garside groups.

\begin{lemma}\label{lem:tinf}
Let $H$ be an abelian subgroup of $G$ which consists of periodic elements.
Then  $\INF|_{H}: H\to \Q$ is a monomorphism.
In particular, $H$ is a cyclic group.
\end{lemma}

\begin{proof}
Let $h_1,h_2\in H$ with $\INF(h_i)=p_i/q_i$ for $i=1,2$.
Because $h_i^{q_i}$ is conjugate to $\Delta^{p_i}$ (by Lemma~\ref{lem:per-elt})
and $\Delta^m$ is central, one has $h_i^{q_im}=\Delta^{p_im}$ for $i=1,2$.
Therefore
$$
(h_1h_2)^{q_1q_2m}
=(h_1^{q_1m})^{q_2} \cdot (h_2^{q_2m})^{q_1}
=\Delta^{p_1mq_2}\Delta^{p_2mq_1}
=\Delta^{m(p_1q_2+p_2q_1)},
$$
hence
$\INF(h_1h_2) = m(p_1q_2+p_2q_1) / q_1q_2m = p_1/q_1+p_2/q_2 =\INF(h_1)+\INF(h_2)$.
This means that $\INF|_{H}: H\to \Q$ is a homomorphism.
If $h\in H$ and $\INF(h)=0$, then
$h$ is conjugate to $\Delta^0=1$ by Lemma~\ref{lem:per-elt}, hence
$h=1$.
This means that $\INF|_{H}: H\to \Q$ is injective.

Notice that, for all $g\in G$, $\INF(g)$ is rational of the form $p/q$ with $1\le q\le\Vert\Delta\Vert$
(see Proposition~\ref{prop:BasicTinf}).
Therefore $\INF(H)$ is a discrete subgroup of $\Q$, hence it is a cyclic group.
Because  $\INF|_{H}: H\to \Q$ is injective,
$H$ is also a cyclic group.
\end{proof}

\begin{theorem}\label{thm:cyclic}
Let $G$ be a Garside group with Garside element $\Delta$.
Then every finite subgroup of\/ $G_\Delta$ is cyclic.
\end{theorem}

\begin{proof}
Let $K$ be a finite subgroup of $G_\Delta$.
Let $H$ be the preimage of $K$ under the natural projection $G\to G_\Delta$.
Notice that every element of $H$ is periodic
and that $H$ is abelian by Theorem~\ref{thm:CMW}.
By Lemma~\ref{lem:tinf}, $H$ is a cyclic group,
hence $K$ is cyclic.
\end{proof}

We give the following for later use.

\begin{lemma}\label{lem:BCMW-a}
Let $g$ be a nonidentity $p/q$-periodic element of $G$.
\begin{enumerate}
\item
$\bar g$ has order $\frac{qm}{p\wedge m}$ in $G_\Delta$.

\item
$\langle \bar g\rangle=\langle\bar g^r\rangle$ in $G_\Delta$
if and only if\/
$r$ is coprime to $\frac{qm}{p\wedge m}$.

\item
Suppose that $\langle \bar g\rangle=\langle\bar g^r\rangle$ in $G_\Delta$
for a nonzero integer $r$.
For $h,x\in G$,
$h=x^{-1}gx$ if and only if $h^r=x^{-1}g^rx$.
Therefore, the CDP and CSP for $(g,h)$
are equivalent to those for $(g^r,h^r)$,
and the centralizer of\/ $g$ in $G$ is the same as
that of\/ $g^r$ in $G$.
\end{enumerate}
\end{lemma}

\begin{proof}
Since $g$ is not the identity, $p$ is not zero.

(i)\ \
By Lemma~\ref{lem:PerCyclic}, there is an element $h$ in $G$ with
$\langle \bar g\rangle=\langle \bar h\rangle$ and $\Delta^m = h^{\frac{qm}{p\wedge m}}$.
Therefore $\bar g$ has order $\frac{qm}{p\wedge m}$ in $G_\Delta$.

(ii)\ \
It follows from (i).

(iii)\ \
It is obvious that $x^{-1}gx=h$ implies $x^{-1}g^r x=h^r$.
Conversely, suppose $x^{-1}g^r x=h^r$.
As $\langle\bar g\rangle=\langle\bar g^r\rangle$ in $G_\Delta$,
we have $g^{rk}=\Delta^{m\ell}g$ hence $g^{rk-1}=\Delta^{m\ell}$
for some $k,\ell\in\Z$.
Since $h$ is periodic and
$$
\INF(h^{rk-1})
=(rk-1)/r\cdot \INF(h^r)
=(rk-1)/r\cdot \INF(g^r)
=\INF(g^{rk-1})=m\ell,
$$
$h^{rk-1}=\Delta^{m\ell}$ by Corollary~\ref{cor:basic},
hence $h^{rk}=\Delta^{m\ell}h$.
As $x^{-1}g^r x=h^r$, we have
$x^{-1}(\Delta^{m\ell}g)x=x^{-1}g^{rk}x=h^{rk}=\Delta^{m\ell}h$.
As $\Delta^{m\ell}$ is central, it follows that $x^{-1}g x=h$.
\end{proof}

\subsection{Slim and precentral}
We define the following notions for periodic elements.

\begin{definition}\label{def:BCMW}
Let $g\in G$ be $p/q$-periodic.
\begin{enumerate}
\item $g$ is said to be \emph{precentral} if $p\equiv 0\bmod m$.
\item $g$ is said to be \emph{slim} if $p\equiv 1\bmod q$.
\end{enumerate}
\end{definition}

If $g\in G$ is $p/q$-periodic then $g^q$ is the minimal
positive power of $g$ which is conjugate to a power of $\Delta$.
Therefore being precentral means this power is central.

\begin{lemma}\label{lem:precentral}
Let $g\in G$ be periodic.
If\/ $g$ is precentral, then so is $g^k$ for all $k\in \Z$.
\end{lemma}

\begin{proof}
Let $g$ be precentral and $p/q$-periodic, then $p\equiv 0\bmod m$.
Choose any $k\in\Z$.
Let $p'=kp/(k\wedge q)$ and
$q'=q/(k\wedge q)$.
Then $\INF(g^k)=\frac{kp}{q}=\frac{kp/(k\wedge q)}{q/(k\wedge q)}=\frac{p'}{q'}$.
Because $(kp)\wedge q=k\wedge q$, $p'$ and $q'$ are coprime, hence
$g^k$ is $p'/q'$-periodic by Corollary~\ref{cor:Periodic_Tinf}.
On the other hand, $p'= p\cdot k/(k\wedge q) \equiv 0\bmod m$ as $p\equiv 0\bmod m$.
Therefore $g^k$ is precentral.
\end{proof}

The terminology `slim' comes from the following observation.

\begin{lemma}\label{lem:Pmin}
Let $g\in G$ be $p/q$-periodic with $q\ge 2$.
Then the following are equivalent.
\begin{enumerate}
\item
$g$ is slim.

\item
$g$ is conjugate to an element of the form $\Delta^u a$
with $\tau^{(q-1)u}(a)\cdots \tau^{u}(a)\,a = \Delta$.

\item
Every element $h\in [g]^{St}$ is of the form $\Delta^u a$
with $\tau^{(q-1)u}(a)\cdots \tau^{u}(a)\,a = \Delta$.
\end{enumerate}
\end{lemma}

\begin{proof}
We prove the equivalences by showing
(iii) $\Rightarrow$  (ii) $\Rightarrow$ (i) $\Rightarrow$ (iii).

(iii) $\Rightarrow$ (ii)\ \
It is obvious.

(ii) $\Rightarrow$ (i)\ \
Since $g$ is $p/q$-periodic, $g^q$ is conjugate to $\Delta^p$.
On the other hand, $g^q$ is conjugate to $(\Delta^u a)^q=\Delta^{uq+1}$ by the hypothesis (ii).
Therefore $\Delta^p = \Delta^{uq+1}$,
hence $p=uq+1\equiv 1\bmod q$.

(i) $\Rightarrow$ (iii)\ \
Note that $p=uq+1$ for some integer $u$.
Then $\SUP(g)=\INF(g)=p/q=u+1/q$ by Lemma~\ref{lem:per-elt}.
Choose any $h\in[g]^{St}$.
Then, by Proposition~\ref{prop:BasicTinf},
$$
\begin{array}{l}
\inf(h)=\infs(g)=\lfloor \INF(g)\rfloor=\lfloor u+1/q\rfloor=u,\\
\sup(h)=\sups(g)=\lceil \SUP(g)\rceil=\lceil u+1/q\rceil=u+1,
\end{array}
$$
from which $h=\Delta^u a$ for some $a\in (1,\Delta)$.
In addition,
$$
\begin{array}{l}
\inf(h^q)=\infs(g^q)=\lfloor q\cdot \INF(g)\rfloor=uq+1,\\
\sup(h^q)=\sups(g^q)=\lceil q\cdot \SUP(g)\rceil=uq+1,
\end{array}
$$
from which $h^q=\Delta^{uq+1}$.
Therefore
$$
\Delta^{uq+1}=h^q=(\Delta^u a)\cdots(\Delta^u a)
=\Delta^{uq}\, \tau^{(q-1)u}(a)\,\tau^{(q-2)u}(a)\cdots \tau^{u}(a)\,a,
$$
which implies $\tau^{(q-1)u}(a)\,\tau^{(q-2)u}(a)\cdots \tau^{u}(a)\,a=\Delta$.
\end{proof}

\begin{lemma}\label{lem:BCMW}
Let $g\in G$ be nonidentity and $p/q$-periodic.
For an integer $r$, $g^r$ is slim with
$\langle \bar g\rangle=\langle\bar g^r\rangle$ in $G_\Delta$
if and only if
$pr\equiv 1\bmod q$ and $r$ is coprime to $\frac{m}{p\wedge m}$.
In particular, such an integer $r$ exists.
\end{lemma}

\begin{proof}
Suppose $g^r$ is slim with
$\langle \bar g\rangle=\langle\bar g^r\rangle$ in $G_\Delta$.
By Lemma~\ref{lem:BCMW-a}(ii),
$r$ is coprime to both $q$ and $\frac{m}{p\wedge m}$.
Note that $\INF(g^r)=\frac{pr}{q}$.
Because $q$ is coprime to both $p$ and $r$,
it is coprime to $pr$,
hence $g^r$ is $pr/q$-periodic.
Therefore $pr\equiv1 \bmod q$ because $g^r$ is slim.

Conversely, suppose $r$ is coprime to $\frac{m}{p\wedge m}$ with
$pr\equiv 1 \bmod q$.
The condition $pr\equiv 1 \bmod q$ implies that $pr$ is coprime to $q$.
Hence $g^r$ is slim because $\INF(g^r)=\frac{pr}{q}$ and
$pr\equiv 1 \bmod q$.
Since $pr$ and $q$ are coprime, so are $r$ and $q$.
Combining with the condition that
$r$ is coprime to $\frac{m}{p\wedge m}$,
we can conclude that $r$ is coprime to $\frac{qm}{p\wedge m}$, hence
$\langle \bar g\rangle=\langle\bar g^r\rangle$
by Lemma~\ref{lem:BCMW-a}(ii).

Since $p$ and $q$ are coprime, so are $p$ and $\frac{qm}{p\wedge m}$.
Therefore there exist integers $r$ and $s$ with $pr+\frac{qm}{p\wedge m}s=1$,
which implies that $pr\equiv 1\bmod q$ and $r$ is coprime to $\frac{m}{p\wedge m}$.
\end{proof}

In the above lemma, the proof shows that it is easy to compute the exponent $r$ of $g^r$,
by applying the Euclidean algorithm to $(p,\frac{qm}{p\wedge m})$.
We remark that Theorem~\ref{thm:CMW}(i) implies
the existence of a slim power $g^r$
with $\langle \bar g\rangle=\langle\bar g^r\rangle$
but does not give the exponent $r$ explicitly.

\subsection{Super summit sets of slim, precentral periodic elements}

Here, we will show that the super summit set of a periodic element has
a useful property provided the element is slim and precentral.
First, we introduce partial cycling, which was extensively
studied by Birman, Gebhardt and Gonz\'alez-Meneses in~\cite{BGG08}.

\begin{definition}\label{defn:pa-cy}
Let $\Delta^u a_1a_2\cdots a_\ell$ be the normal form of $g\in G$.
Let $b\in (1,\Delta)$ be a prefix of $a_1$,
i.e. $a_1=ba_1'$ for some $a'_1\in [1,\Delta)$.
The conjugation
$$
\tau^{-u}(b)^{-1} g \tau^{-u}(b) = \Delta^u a_1'a_2\cdots a_\ell\tau^{-u}(b)
$$
is called the \emph{partial cycling} of $g$ by $b$.
\end{definition}

Any partial cycling does not decrease the infimum,
hence summit sets are closed under any partial cycling.
If $g\in G$ is periodic, then $[g]^S=[g]^U$
because every element in $[g]^S$ has canonical length $\le 1$
by Lemma~\ref{lem:per-elt}.
The following theorem shows a property of slim, precentral periodic elements.

\begin{theorem}\label{thm:pa-cy}
Let $g$ be a slim, precentral periodic element
of a Garside group $G$. Then
$$
[g]^{\inf}=[g]^S=[g]^U=[g]^{St}.
$$
In particular, $[g]^{S}$ is closed under any partial cycling.
\end{theorem}

\begin{proof}
Since $[g]^{\inf}$ is closed under any partial cycling
and $[g]^{St}\subset[g]^S=[g]^U\subset[g]^{\inf}$,
it suffices to show that $[g]^{\inf}\subset [g]^{St}$.

Suppose that $h$ is an element of $[g]^{\inf}$.
Let $g$ be $p/q$-periodic.
If $q=1$, then $g$ is conjugate to $\Delta^p$,
hence there is nothing to prove.
Let $q\ge 2$.
As $g$ is slim and precentral, $p=uq+1=m\ell$
for some integers $u$ and $\ell$.
For all integers $k\ge 1$,
\begin{align*}
\INF(h)&= \SUP(h)=p/q=u+1/q,\\
\infs(h^k)&=  \lfloor k\INF(h)\rfloor
=ku+\lfloor k/q\rfloor,\\
\sups(h^k)&=  \lceil k\SUP(h)\rceil
= ku+\lceil k/q\rceil.
\end{align*}
Because $\INF(h^q)=\INF(g^q)=p$, one has
$$
h^q=\Delta^p=\Delta^{uq+1}=\Delta^{m\ell}
$$
by Corollary~\ref{cor:basic}.
Because $\inf(h)=\infs(h)=u$, there exists a positive element $a$ such that
$$
h=\Delta^u a.
$$
Because $\sup(h)\ge \sups(h)=u+1$, $a$ is not the identity.
Let $\psi$ denote $\tau^{u}$. For all integers $k\ge 1$,
$$
h^k=\Delta^{ku}\psi^{k-1}(a)\psi^{k-2}(a)\cdots\psi(a) a.
$$
Since $h^q=\Delta^{qu+1}$, one has
$$
\Delta^{qu+1}=h^q=\Delta^{qu}\psi^{q-1}(a)\psi^{q-2}(a)\cdots\psi(a) a.
$$
Hence $\Delta=\psi^{q-1}(a)\psi^{q-2}(a)\cdots\psi(a) a$.
In particular,
$\psi^{k-1}(a)\psi^{k-2}(a)\cdots a \in (1,\Delta)$
for all integers $k$ with $1\le k < q$.
Therefore
$$
\inf(h^k)=uk=\infs(h^k)
\quad\mbox{and}\quad
\sup(h^k)=uk+1=\sups(h^k)
\qquad\mbox{for $k=1,\ldots,q-1$}.
$$
Because $h^q=\Delta^{uq+1}$, this proves that
$h\in[g]^{St}$.
\end{proof}

The following example shows that both slimness and precentrality
are indeed necessary in the above theorem.

\begin{example}\label{eg:conditions}
Let the $n$-braid group $B_n$ be endowed with the classical Garside structure.
We will write $\varepsilon=\varepsilon_{(n)}$ in order to specify the braid index.
Since $\varepsilon_{(n)}^{n-1}=\Delta^2$,
$\INF(\varepsilon_{(n)})=2/(n-1)$.
Consider
$\varepsilon_{(5)}=(\sigma_4\sigma_3\sigma_2\sigma_1)\sigma_1\in B_5$
and $\varepsilon_{(6)}=(\sigma_5\sigma_4\sigma_3\sigma_2\sigma_1)\sigma_1\in B_6$.
Since
$$
\INF(\varepsilon_{(5)}) = 2/4=1/2
\quad\mbox{and}\quad
\INF(\varepsilon_{(6)}) = 2/5,
$$
$\varepsilon_{(5)}$ is slim but not precentral,
whereas $\varepsilon_{(6)}$ is precentral but not slim.
In addition,
\begin{align*}
& \inf(\varepsilon_{(5)})=\inf(\varepsilon_{(6)})=0
\quad\mbox{and}\quad
\len(\varepsilon_{(5)})=\len(\varepsilon_{(6)})=2;\\
& \infs(\varepsilon_{(5)})=\infs(\varepsilon_{(6)})=0
\quad\mbox{and}\quad
\lens(\varepsilon_{(5)})=\lens(\varepsilon_{(6)})=1;\\
& \lens(\varepsilon_{(5)}^2)=0
\quad\mbox{and}\quad
\lens(\varepsilon_{(6)}^2)=1.
\end{align*}

From the above identities,
$\varepsilon_{(5)}\in[\varepsilon_{(5)}]^{\inf}\setminus[\varepsilon_{(5)}]^S$,
hence $[\varepsilon_{(5)}]^{\inf}\ne[\varepsilon_{(5)}]^S$.
Similarly, $[\varepsilon_{(6)}]^{\inf}\ne[\varepsilon_{(6)}]^S$.

Consider the elements
$g_1=\sigma_1(\sigma_4\sigma_3\sigma_2\sigma_1)\in B_5$
and $g_2=\sigma_1(\sigma_5\sigma_4\sigma_3\sigma_2\sigma_1)\in B_6$.
The partial cycling on $g_1$ and $g_2$ by $\sigma_1$ yields
$\sigma_1^{-1}g_1\sigma_1=\varepsilon_{(5)}$ and
$\sigma_1^{-1}g_2\sigma_1=\varepsilon_{(6)}$, respectively.
Since $\inf(g_1)=\inf(g_2)=0$ and $\len(g_1)=\len(g_2)=1$,
we have $g_1\in [\varepsilon_{(5)}]^S$ and $g_2\in [\varepsilon_{(6)}]^S$.
Because neither $\varepsilon_{(5)}$ nor $\varepsilon_{(6)}$ is a super summit element,
this shows that neither $[\varepsilon_{(5)}]^S$ nor $[\varepsilon_{(6)}]^S$
is closed under partial cycling.

The normal forms of $g_1^2$ and $g_2^2$ are as in the right hand sides
in the following equations:
\begin{align*}
g_1^2
&=(\sigma_1\sigma_4\sigma_3\sigma_2\sigma_1)
   (\sigma_1\sigma_4\sigma_3\sigma_2\sigma_1)
  =(\sigma_1\sigma_4\sigma_3\sigma_2\sigma_1\sigma_4\sigma_3\sigma_2)
   \cdot(\sigma_1\sigma_2),\\
g_2^2
&=(\sigma_1\sigma_5\sigma_4\sigma_3\sigma_2\sigma_1)
   (\sigma_1\sigma_5\sigma_4\sigma_3\sigma_2\sigma_1)
  =(\sigma_1\sigma_5\sigma_4\sigma_3\sigma_2\sigma_1\sigma_5\sigma_4\sigma_3\sigma_2)
   \cdot(\sigma_1\sigma_2).
\end{align*}
In particular, both $g_1^2$ and $g_2^2$ have canonical length 2, hence
they do not belong to their super summit sets.
This means that $g_1$ and $g_2$ do not belong to their stable super summit sets.
Hence $[\varepsilon_{(5)}]^S\ne[\varepsilon_{(5)}]^{St}$
and $[\varepsilon_{(6)}]^S\ne[\varepsilon_{(6)}]^{St}$.
\end{example}

\smallskip

Before closing this section, we make some remarks on
the requirement of being slim and precentral in Theorem~\ref{thm:pa-cy}.

Given a periodic element $g\in G$,
it is easy to compute a nonzero integer $r$
such that $g^r$ is slim
and for $h,x\in G$, $h=x^{-1}gx$ if and only if $h^r=x^{-1}g^rx$,
by Lemmas~\ref{lem:BCMW-a} and~\ref{lem:BCMW}.
Taking such a power of $g$,
we may assume without loss of generality that $g$ is slim
when thinking about the conjugacy problem for $g$.

For the precentrality condition,
we can make every periodic element precentral
by modifying the Garside structure on $G$:
if $(G^+,\Delta)$ is a Garside structure on $G$, then
$(G^+,\Delta^m)$ is also a Garside structure on $G$
under which every periodic element is precentral.

\def\temp{For the precentrality condition,
we can make every periodic element precentral
by modifying the Garside structure on $G$:
for a Garside group $G$ with Garside element $\Delta$,
if we change the Garside structure on $G$ by declaring
a central power $\Delta^m$ as a new Garside element, then every periodic element
becomes precentral in this new Garside structure.
}

\section{Periodic elements in some Garside groups arising from reflection groups}
\label{sec:SomeGroups}

This section studies periodic elements in
the Artin groups of type $\arA_n$, $\arB_n$, $\arD_n$, $\arI_2(e)$
and the braid group of the complex reflection group of type $(e,e,n)$.
These groups are known to be Garside groups.
Using the recent result of Bessis~\cite{Bes06a} on
periodic elements in the braid groups of complex reflection groups,
we will find all the primitive periodic elements, and then investigate precentrality
for periodic elements in those Garside groups.

First we review Artin groups and braid groups of complex reflection groups.
See~\cite{BS72, Del72, Hum90} for Artin groups and
\cite{Bro01} for braid groups of complex reflection groups.
Let $\langle ab\rangle^k$ denote the alternating product
$abab\cdots$ of length $k$.
For instance, $\langle ab\rangle^3 = aba$.

\subsection{Artin groups}
Let $M$ be a symmetric $n\times n$ matrix with entries $m_{ij}\in\Z_{\ge 1}\cup\{\infty\}$
where $m_{ii}=1$ and $m_{ij}\ge 2$ for $i\ne j$.
The \emph{Artin group} of type $M$ is defined by the presentation
\begin{equation}\label{eq:Artin}
A(M)=\langle s_1,\ldots,s_n\mid \langle s_is_j\rangle^{m_{ij}}
=\langle s_js_i\rangle^{m_{ji}}
\quad\mbox{for all $i\ne j$ with $m_{ij}\ne\infty$}\rangle.
\end{equation}

The \emph{Coxeter group} $W(M)$ of type $M$ is the quotient
of $A(M)$ by the relations $s_i^2=1$.
We say that the Artin group $A(M)$ is \emph{of finite type} if
the associated Coxeter group $W(M)$ is a finite set.

\begin{figure}
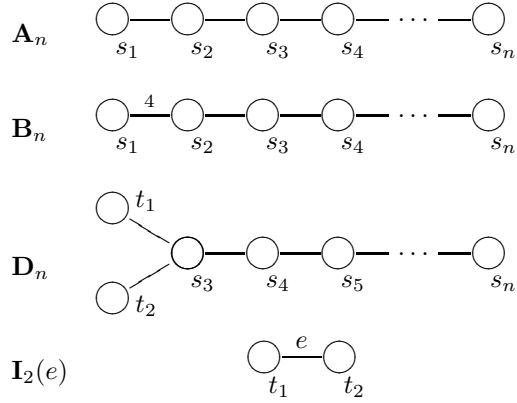
\def\arraystretch{1.5}
\begin{tabular}{lc}
$\arA_n$ &
$\xy
( 0,2)  *++={\phantom{2}} *\frm{o} **@{-};
(10,2)  *++={\phantom{2}} *\frm{o} **@{-};
(20,2)  *++={\phantom{2}} *\frm{o} **@{-};
(30,2)  *++={\phantom{2}} *\frm{o} **@{-};
(40,2)  *++={\dots} **@{-} ;
(50,2)  *++={\phantom{2}} *\frm{o} **@{-};
( 2,-2) *++={s_1};
(12,-2) *++={s_2};
(22,-2) *++={s_3};
(32,-2) *++={s_4};
(52,-2) *++={s_n};
\endxy$\\[1em]
$\arB_n$ &
$\xy
( 0,2)  *++={\phantom{2}} *\frm{o} **@{-};
(10,2)  *++={\phantom{2}} *\frm{o} **@{-};
(20,2)  *++={\phantom{2}} *\frm{o} **@{-};
(30,2)  *++={\phantom{2}} *\frm{o} **@{-};
(40,2)  *++={\dots} **@{-} ;
(50,2)  *++={\phantom{2}} *\frm{o} **@{-};
( 5, 4) *++={{}_4};
( 2,-2) *++={s_1};
(12,-2) *++={s_2};
(22,-2) *++={s_3};
(32,-2) *++={s_4};
(52,-2) *++={s_n};
\endxy$\\[1em]
$\arD_n$ &
$\xy
(10,8)  *++={\phantom{2}} *\frm{o} ;
(20,2)  *++={\phantom{2}} *\frm{o} **@{-};
(10,-4) *++={\phantom{2}} *\frm{o} ;
(20,2)  *++={\phantom{2}} *\frm{o} **@{-};
(30,2)  *++={\phantom{2}} *\frm{o} **@{-};
(40,2)  *++={\phantom{2}} *\frm{o} **@{-};
(50,2)  *++={\dots} **@{-} ;
(60,2)  *++={\phantom{2}} *\frm{o} **@{-};
(14.5, 9)   *++={t_1} ;
(14.5,-5)  *++={t_2} ;
(22,-2) *++={s_3};
(32,-2) *++={s_4};
(42,-2) *++={s_5};
(62,-2) *++={s_n};
\endxy$\\[1em]
$\arI_2(e)$ &
$\xy
( 0,2)  *++={\phantom{2}} *\frm{o} **@{-};
(10,2)  *++={\phantom{2}} *\frm{o} **@{-};
( 5, 4) *++={e};
( 2,-2) *++={t_1};
(12,-2) *++={t_2};
\endxy$
\end{tabular}
\caption{Coxeter graphs}
\label{fig:graph}
\end{figure}

It is convenient to define an Artin group by a \emph{Coxeter graph},
whose vertices are labeled by the generators $s_1,\ldots, s_n$
and which has an edge labeled $m_{ij}$
between the vertices $s_i$ and $s_j$ whenever $m_{ij}\ge 3$
or $m_{ij}=\infty$.
The label 3 is usually suppressed.
The Coxeter graphs of type $\arA_n$, $\arB_n$, $\arD_n$ and $\arI_2(e)$
are in Figure~\ref{fig:graph}.
The Artin groups of these types are of finite type,
and the associated Coxeter groups are real reflection groups.
We denote the generators of these Artin groups as in Figure~\ref{fig:graph}.

Since the relations in~(\ref{eq:Artin}) involve only positive words,
the presentation defines a monoid.
We denote this monoid by $A(M)^+$.
In other words, $A(M)^+$ consists of positive words in the generators
modulo the defining relations.

By the study of Brieskorn and Saito~\cite{BS72} and Deligne~\cite{Del72},
it is known that if an Artin group $A(M)$ is of finite type, then
it is a Garside group with Garside monoid $A(M)^+$.
The Garside element $\Delta$ is the least common multiple of the generators
in the presentation~(\ref{eq:Artin}).
For instance,
$$
\begin{array}{ll}
\Delta=s_1(s_2s_1)\cdots(s_ns_{n-1}\cdots s_1)&\mbox{in $A(\arA_n)$;}\\
\Delta=(s_ns_{n-1}\cdots s_1)^n               &\mbox{in $A(\arB_n)$;}\\
\Delta=(s_n s_{n-1}\cdots s_3 t_1t_2)^{n-1}    &\mbox{in $A(\arD_n)$;}\\
\Delta=\langle t_1t_2\rangle^e         &\mbox{in $A(\arI_2(e))$.}
\end{array}
$$
We refer to this Garside structure $(A(M)^+, \Delta)$
as the \emph{classical Garside structure} on $A(M)$.

Artin groups of finite type have another Garside structure,
called the \emph{dual Garside structure}.
This structure was constructed originally by Birman, Ko and Lee~\cite{BKL98}
for $A(\arA_n)$,
and then by Bessis~\cite{Bes03} for all finite type Artin groups.
In the construction of Bessis, a choice of a Coxeter element
in the associated Coxeter group determines the dual Garside structure,
in particular, the Garside element $\delta$.
We choose the Garside element as follows:
$$
\begin{array}{ll}
\delta=s_ns_{n-1}\cdots s_1     &\mbox{in $A(\arA_n)$ and $A(\arB_n)$;}\\
\delta=s_ns_{n-1}\cdots s_3t_1t_2     &\mbox{in $A(\arD_n)$;}\\
\delta=t_1t_2                   &\mbox{in $A(\arI_2(e))$.}
\end{array}
$$

From now on, we assume $n\ge 2$ for $A(\arA_n)$ and $A(\arB_n)$,
$n\ge 3$ for $A(\arD_n)$ and $e\ge 3$ for $A(\arI_2(e))$.

\subsection{Braid groups of complex reflection groups}

Let $V$ be a finite dimensional complex vector space.
A \emph{complex reflection group} in $GL(V)$ is a subgroup $W$
of the general linear group $GL(V)$
generated by complex reflections---nontrivial elements of $GL(V)$ that fix a complex
hyperplane in $V$ pointwise.
Irreducible finite complex reflection groups were classified
by Shephard and Todd~\cite{ST54}.
There are a general infinite family $G(de,e,n)$
for $d,e,n\in\Z_{\ge 1}$, and 34 exceptions labeled $G_4,\ldots, G_{37}$.
See~\cite{BMR98, Bro01, BM04, Bes06a} for the presentations
of the complex reflection groups and their braid groups.
Special cases of complex reflection groups are
isomorphic to real reflection groups:
\begin{itemize}
\item[]
$G(1,1,n)$ is the Coxeter group of type $\arA_{n-1}$;\\
$G(2,1,n)$ is the Coxeter group of type  $\arB_n$;\\
$G(2,2,n)$ is the Coxeter group of type  $\arD_n$;\\
$G(e,e,2)$ is the Coxeter group of type  $\arI_2(e)$.
\end{itemize}

\smallskip
Let $V'$ be the complement of all reflecting hyperplanes of
reflections in a complex reflection group $W\subset GL(V)$.
Then $W$ acts on $V'$.
The fundamental group $\pi_1(W\backslash V')$
of the quotient space $W\backslash V'$
is called the \emph{braid group} of $W$, denoted $B(W)$.
The fundamental group $\pi_1(V')$ is called the
\emph{pure braid group} of $W$, denoted $P(W)$.
Let $B(de,e,n)$ denote the braid group of the complex refection group $G(de,e,n)$.

This paper is interested in the braid groups $B(e,e,n)$ for $e\ge 2$ and $n\ge 3$,
which are known to be Garside groups.
The Artin group $A(\arD_n)$ is a special case of $B(e,e,n)$ with $e=2$.
It is not known whether the braid groups $B(de,e,n)$ with $d\ge 2$ have
a Garside structure.

Bessis and Corran~\cite{BC06} constructed
the dual Garside structure for $B(e,e,n)$,
and then Bessis~\cite{Bes06a} improved this result,
giving a new geometric interpretation and extending the construction
to the exceptional cases not covered before.
Bessis and Corran also showed in~\cite{BC06} that the monoid arising from the
Brou\'e-Malle-Rouquier presentation for $B(e,e,n)$ in~\cite{BMR98}
is not a Garside monoid.
Combining the Brou\'e-Malle-Rouquier presentation and the Bessis-Corran presentation,
Corran and Picantin~\cite{CP09} recently proposed a presentation for $B(e,e,n)$
which gives a new Garside structure.

In the Brou\'e-Malle-Rouquier presentation,
$B(e,e,n)$ is generated by
$t_1,t_2,s_3,s_4,\ldots,s_n$ with the following defining relations:
\begin{itemize}
\item[] $s_is_j=s_js_i$ if $|i-j|\ge 2$;
\item[] $s_is_{i+1}s_i=s_{i+1}s_is_{i+1}$ for $i=3,\ldots,n-1$;
\item[] $\langle t_1t_2\rangle^e = \langle t_2t_1\rangle^e$;
\item[] $t_1s_j=s_jt_1$ and $t_2s_j=s_jt_2$ for $j\ge 4$;
\item[] $t_1s_3t_1=s_3t_1s_3$ and $t_2s_3t_2=s_3t_2s_3$;
\item[] $s_3t_1t_2s_3t_1t_2=t_1t_2s_3t_1t_2s_3$.
\end{itemize}

\begin{figure}
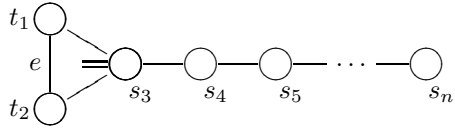

$$
\xy
(10,8)  *++={\phantom{2}} *\frm{o};
(10,-4) *++={\phantom{2}} *\frm{o} **@{-};
(13,2)  *++={} ;
(20,2)  *++={\phantom{2}} *\frm{o} **@{=};
(10,8)  *++={\phantom{2}} *\frm{o} ;
(20,2)  *++={\phantom{2}} *\frm{o} **@{-};
(10,-4) *++={\phantom{2}} *\frm{o} ;
(20,2)  *++={\phantom{2}} *\frm{o} **@{-};
(30,2)  *++={\phantom{2}} *\frm{o} **@{-};
(40,2)  *++={\phantom{2}} *\frm{o} **@{-};
(50,2)  *++={\dots} **@{-} ;
(60,2)  *++={\phantom{2}} *\frm{o} **@{-};
(6,8)   *++={t_1} ;
(6,-4)  *++={t_2} ;
(22,-2) *++={s_3};
(32,-2) *++={s_4};
(42,-2) *++={s_5};
(62,-2) *++={s_n};
(8,2)   *++={e}
\endxy
$$
\caption{Brou\'e-Malle-Rouquier presentation for $B(e,e,n)$}
\label{fig:BMR-rel}
\end{figure}

This presentation is usually illustrated as in Figure~\ref{fig:BMR-rel},
which looks like a Coxeter graph.
In the figure, the symbol
``\raise1pt\hbox to 0pt{\rule{9pt}{.5pt}\hss}%
\raise3pt\hbox{\rule{9pt}{.5pt}}''
at the vertex $s_3$ indicates the relation
$s_3t_1t_2s_3t_1t_2= t_1t_2s_3t_1t_2s_3$.

\smallskip
In the Corran-Picantin presentation,
$B(e,e,n)$ is generated by $t_1,\ldots,t_e,s_3,\ldots,s_n$
with the following defining relations:
\begin{itemize}
\item[] $s_is_j=s_js_i$ if $|i-j|\ge 2$;
\item[] $s_is_{i+1}s_i=s_{i+1}s_is_{i+1}$ for $i=3,\ldots,n-1$;
\item[] $t_1t_2=t_2t_3=\cdots=t_{e-1}t_e=t_et_1$;
\item[] $t_is_j=s_jt_i$ for $i=1,\ldots,e$ and $j\ge 4$;
\item[] $t_is_3t_i = s_3 t_i s_3$ for $i=1,\ldots,e$.
\end{itemize}

\begin{figure}
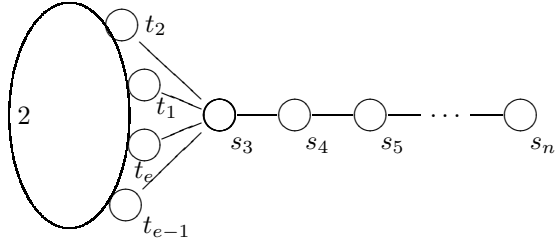

$$\xy
{\ellipse(8,15){}};
(7, 12) *++={\phantom{2}} *\frm{o};
(20,0) *++={\phantom{2}} *\frm{o} **@{-};
(7.5,-12) *++={\phantom{2}} *\frm{o};
(20,0) *++={\phantom{2}} *\frm{o} **@{-};
(10, 4) *++={\phantom{2}} *\frm{o};
(20,0) *++={\phantom{2}} *\frm{o} **@{-};
(10,-4) *++={\phantom{2}} *\frm{o};
(20,0) *++={\phantom{2}} *\frm{o} **@{-};
(30, 0) *++={\phantom{2}} *\frm{o} **@{-};
(40, 0) *++={\phantom{2}} *\frm{o} **@{-};
(50, 0) *++={\dots}  **@{-};
(60, 0) *++={\phantom{2}} *\frm{o} **@{-};
(-6,0) *++={2};
(11.5, 12) *++={t_2};
(13,-15) *++={t_{e-1}};
(13, 1) *++={t_1};
(10,-7.5) *++={t_e};
(23,-4) *++={s_3};
(33,-4) *++={s_4};
(43,-4) *++={s_5};
(63,-4) *++={s_n}
\endxy
$$
\caption{Corran-Picantin presentation for $B(e,e,n)$}\label{fig:New-rel}
\end{figure}

This presentation is illustrated in Figure~\ref{fig:New-rel},
where the large circle with label 2
indicates the relation $t_1t_2=t_2t_3=\cdots=t_{e-1}t_e=t_et_1$.
We remark that the above presentation is slightly different from
but equivalent to the one given by Corran and Picantin:
in their presentation the generators $t_1,\ldots,t_e$ satisfy
$t_2t_1=t_3t_2=\cdots=t_et_{e-1}=t_1t_e$ rather than
$t_1t_2=t_2t_3=\cdots=t_{e-1}t_e=t_et_1$.

One can understand the difference between the Brou\'e-Malle-Rouquier
and Corran-Picantin presentations as follows.
The subgraph
$\small\xy
(20,2)  *++={\phantom{2}} *\frm{o} ;
(30,2)  *++={\phantom{2}} *\frm{o} **@{-};
(25,-1) *++={e};
(20,-2) *++={t_1};
(30,-2) *++={t_2}
\endxy$
in Figure~\ref{fig:BMR-rel} gives a presentation of the Artin group $A(\arI_2(e))$,
which gives the classical Garside structure on $A(\arI_2(e))$.
In Figure~\ref{fig:New-rel}, this part is replaced by a subgraph which gives
the dual Garside structure on $A(\arI_2(e))$.
Therefore the Corran-Picantin presentation may be regarded as
a mixture of the Brou\'e-Malle-Rouquier presentation
and the Bessis-Corran dual presentation.
Compared to the dual Garside structure,
we will refer to the Garside structure arising from the Corran-Picantin presentation
as the \emph{classical Garside structure}.

\smallskip
Let $S$ denote the word $s_ns_{n-1}\cdots s_3$.
The Garside element $\Delta$ in the classical Garside structure
and  the Garside element $\delta$
in the dual Garside structure are as follows:
$$
\Delta=(S t_1t_2)^{n-1}
\quad\mbox{and}\quad
\delta=S t_1t_2.
$$

From now on, we assume $e\ge 2$ and $n\ge 3$ for $B(e,e,n)$.

\subsection{Result of Bessis}

Let $W\subset GL(V)$ be a complex reflection group.
The largest degree is called the \emph{Coxeter number} of $W$,
which we will denote by $h$.
An integer $d$ is called a \emph{regular number}
if there exist an element $w\in W$ and a complex $d$-th root $\zeta$
of unity such that $\ker(w-\zeta)\cap V'\ne\emptyset$,
where $V'$ is the complement of all reflecting hyperplanes of
reflections in $W$.

\smallskip
The following theorem collects Bessis' results which we will use;
see Lemma 6.11 and Theorems 1.9, 8.2, 12.3, 12.5 in \cite{Bes06a}.

\begin{theorem}[\cite{Bes06a}]\label{thm:bes06a}
Let\/ $W$ be an irreducible well-generated complex reflection group,
with degrees $d_1,\ldots,d_n$, codegrees $d_1^*,\ldots,d_n^*$
and Coxeter number $h$.
Then its braid group $B(W)$ admits the dual Garside structure with
Garside element $\delta$, and the following hold.

\begin{enumerate}
\item
The element $\mu=\delta^h$ is central in $B(W)$ and lies in the pure braid group $P(W)$.

\item
Let $h'=h/(d_1\wedge\cdots\wedge d_n)$.
The center of $B(W)$ is a cyclic group generated by $\delta^{h'}$.

\item
Let $d$ be a positive integer, and let
$$
A(d)=\{1\le i\le n: d|d_i\}
\quad\mbox{and}\quad
B(d)=\{1\le i\le n: d|d_i^*\}.
$$
Then $|A(d)|\le|B(d)|$, and the following conditions are equivalent:
\begin{enumerate}
\item
$|A(d)|=|B(d)|$;

\item
there exists a $d$-th root of $\mu$;

\item
$d$ is regular.
\end{enumerate}
Moreover, the $d$-th root of $\mu$, if exists, is unique up to conjugacy in $B(W)$.
\end{enumerate}
\end{theorem}

The groups $W(\arA_n), W(\arB_n)$ and $G(e,e,n)$ including $W(\arD_n)$ and $W(\arI_2(e))$
are all irreducible well-generated complex reflection groups.

\subsection{Primitive periodic elements}

In \S\ref{ssec:Primitive}, we have shown that every primitive periodic element
in a Garside group is
a root of $\Delta^m$, where $\Delta$ is the Garside element and $\Delta^m$
is the minimal positive power of $\Delta$ which is central.
In this subsection, we give explicitly all the primitive periodic elements in the groups
$A(\arA_n)$, $A(\arB_n)$, $A(\arD_n)$, $A(\arI_2(e))$ and $B(e,e,n)$.
These groups have cyclic centers,
hence the periodicity of an element
does not depend on the choice of a particular Garside structure.

\begin{lemma}\label{lem:reg-num}
Let $d$ be a positive integer.
\begin{enumerate}
\item
In $A(\arA_n)$, $d$ is regular
if and only if\/ $d| n$ or $d| (n+1)$.

\item
In $A(\arB_n)$, $d$ is regular
if and only if\/ $d| (2n)$.

\item
In $A(\arD_n)$, $d$ is regular
if and only if\/ $d| n$ or $d| 2(n-1)$.

\item
In $A(\arI_2(e))$, $d$ is regular
if and only if\/ $d| 2$ or $d| e$.

\item
In $B(e,e,n)$, $d$ is regular
if and only if\/ $d| n$ or $d| e(n-1)$.
\end{enumerate}
\end{lemma}

\begin{proof}
It is known that for well-generated complex reflection groups
the codegrees $d_1^*, \ldots, d_n^*$ are related to the degrees $d_1, \ldots, d_n$
by the formula $d_i+d_i^*=d_n$ for all $i$.
In $A(\arA_n)$, the following are known.
\begin{align*}
\{d_1,\ldots,d_n\}    &= \{2,3,\ldots,n-1\}\cup\{n,n+1\}\\
\{d_1^*,\ldots,d_n^*\}&= \{2,3,\ldots,n-1\}\cup\{0,1\}
\end{align*}
Thus, an integer $d$ is regular if and only if two sets $\{n,n+1\}$
and $\{0,1\}$ have the same number of multiples of $d$,
and this happens if and only if either $d| n$ or $d| (n+1)$.
This proves (i).

The other statements (ii)--(v) can be proved similarly
by using the degrees and codegrees in Table~\ref{tab:deg}.
See~\cite{Hum90} for the degrees of Coxeter groups
and~\cite{Bro01} for the degrees and codegrees of the complex reflection
group $G(e,e,n)$.
\end{proof}

\begin{table}
$$\arraycolsep=3pt
\begin{array}{c|l|c|c}\thickhline
\mbox{Groups}
  & \hspace{7em}\mbox{\centering Degrees and codegrees}
  & d
  & h
  \\\thickhline
\begin{array}{c}
A(\arA_n) \\
\mbox{\scriptsize $(n\ge 2)$}
\end{array}
  & \arraycolsep=0pt\begin{array}{l}
    \{d_1,\ldots,d_n\}    = \{2,3,\ldots,n-1\}\cup\{n,n+1\}\\
    \{d_1^*,\ldots,d_n^*\}= \{2,3,\ldots,n-1\}\cup\{0,1\}
    \end{array}
  & \arraycolsep=0pt\begin{array}{c}
    d| n\\
    d| (n+1)
    \end{array}
  & n+1
  \\\hline
\begin{array}{c}
A(\arB_n) \\
\mbox{\scriptsize $(n\ge 2)$}
\end{array}
  & \arraycolsep=0pt\begin{array}{l}
    \{d_1,\ldots,d_n\}    = \{2,4,\ldots,2n-2\}\cup\{2n\}\\
    \{d_1^*,\ldots,d_n^*\}= \{2,4,\ldots,2n-2\}\cup\{0\}
    \end{array}
  & d| (2n)
  & 2n
  \\\hline
\begin{array}{c}
A(\arD_n) \\
\mbox{\scriptsize $(n\ge 3)$}
\end{array}
  & \arraycolsep=0pt\begin{array}{l}
    \{d_1,\ldots,d_n\}    = \{2,4,\ldots,2n-4\}\cup\{2n-2,n\}\\
    \{d_1^*,\ldots,d_n^*\}= \{2,4,\ldots,2n-4\}\cup\{0,n-2\}
    \end{array}
  & \arraycolsep=0pt\begin{array}{c}
    d| n\\
    d| 2(n-1)
    \end{array}
  & 2(n-1)
  \\\hline
\begin{array}{c}
A(\arI_2(e)) \\
\mbox{\scriptsize $(e\ge 3)$}
\end{array}
  & \arraycolsep=0pt\begin{array}{l}
    \{d_1, d_2\}     = \{2,e\}\\
    \{d_1^*, d_2^*\} = \{0,e-2\}
    \end{array}
  & \arraycolsep=0pt\begin{array}{c}
    d| 2\\
    d| e
    \end{array}
  & e \\\hline
\begin{array}{c}
B(e,e,n) \\
\mbox{\scriptsize $(e\ge 2, n\ge 3)$}
\end{array}
  & \arraycolsep=0pt\begin{array}{l}
    \{d_1,\ldots,d_n\}    = \{e,2e,\ldots,(n-2)e\}\cup\{(n-1)e,n\}\\
    \{d_1^*,\ldots,d_n^*\}= \{e,2e,\ldots,(n-2)e\}\cup\{0,(n-1)e-n\}
    \end{array}
  & \arraycolsep=0pt\begin{array}{c}
    d| n\\
    d| e(n-1)
    \end{array}
  & e(n-1) \\\thickhline
\end{array}
$$
\caption{Degrees, codegrees, regular numbers $d$ and Coxeter number $h$}
\label{tab:deg}
\end{table}

Recall the periodic braid
$\varepsilon = (\sigma_{n-1}\cdots\sigma_1)\sigma_1$ in the braid group $B_n$.
For the groups $A(\arD_n)$ and $B(e,e,n)$,
we define elements $\varepsilon$ as follows:
$$
\begin{array}{ll}
\varepsilon=St_1 S t_2                      & \mbox{in $A(\arD_n)$;}\\
\varepsilon=S t_1 S t_2 \cdots S t_e        & \mbox{in $B(e,e,n)$,}
\end{array}
$$
where $S$ denotes the word $s_ns_{n-1}\cdots s_3$.

The following lemma shows that the elements $\varepsilon$ are also periodic
in both $A(\arD_n)$ and $B(e,e,n)$ like in $B_n$.
The definitions of $\delta$, $\varepsilon$ and $\Delta$ together with
some relations between them are collected in Table~\ref{tab:rel}.

\begin{table}
$$\arraycolsep=3pt
\begin{array}{c|l|c|c|c}\thickhline
\mbox{Groups}
& \hspace{2em}\mbox{periodic elements}
& h'
& \delta^h
& \delta^{h'}
\\\thickhline
\begin{array}{c}
A(\arA_n) \\
\mbox{\scriptsize $(n\ge 2)$}
\end{array}
&
  \begin{array}{l}
  \Delta=s_1(s_2s_1)\cdots(s_n\cdots s_1)\\
  \delta=s_ns_{n-1}\cdots s_1\\
  \varepsilon=(s_ns_{n-1}\cdots s_1)s_1
  \end{array}
    & n+1 &
  \delta^{n+1}=\Delta^2=\varepsilon^n &
  \delta^{n+1}=\Delta^2=\varepsilon^n \\\hline
\begin{array}{c}
A(\arB_n) \\
\mbox{\scriptsize $(n\ge 2)$}
\end{array}
&
  \begin{array}{l}
    \Delta=(s_ns_{n-1}\cdots s_1)^n\\
    \delta=s_ns_{n-1}\cdots s_1
    \end{array}
  & n   &
  \delta^{2n}=\Delta^2 &
  \delta^n=\Delta \\\hline
\begin{array}{c}
A(\arD_n) \\
\mbox{\scriptsize $(n\ge 3)$}
\end{array}
&
  \begin{array}{l}
    \Delta=(St_1t_2)^{n-1}\\
    \delta=St_1t_2\\
    \varepsilon=St_1St_2
    \end{array}
  & \frac{2(n-1)}{2\wedge n} &
  \delta^{2(n-1)} = \Delta^2 = \varepsilon^n &
  \delta^{\frac{2(n-1)}{2\wedge n}}=\Delta^{\frac{2}{2\wedge n}}
  = \varepsilon^{\frac{n}{2\wedge n}} \\\hline
\begin{array}{c}
A(\arI_2(e)) \\
\mbox{\scriptsize $(e\ge 3)$}
\end{array}
&
  \begin{array}{l}
    \Delta=\langle t_1t_2\rangle^e\\
    \delta=t_1t_2\\
    \end{array}
  & \frac{e}{e\wedge 2} &
  \delta^e =\Delta^2 &
  \delta^{\frac{e}{e\wedge 2}} =\Delta^{\frac{2}{e\wedge 2}}\\\hline
\begin{array}{c}
B(e,e,n) \\
\mbox{\scriptsize $(e\ge 2, n\ge 3)$}
\end{array}
&
  \begin{array}{l}
    \Delta=(St_1t_2)^{n-1}\\
    \delta=St_1t_2\\
    \varepsilon=St_1St_2\cdots St_e
    \end{array}
  & \frac{e(n-1)}{e\wedge n} &
  \delta^{e(n-1)}=\Delta^e=\varepsilon^n &
  \delta^{\frac{e(n-1)}{e\wedge n}}
  =\Delta^{\frac{e}{e\wedge n}}
  =\varepsilon^{\frac{n}{e\wedge n}}\\\thickhline
\end{array}
$$
\caption{Periodic elements $\delta$, $\varepsilon$ and $\Delta$, where
$S=s_ns_{n-1}\cdots s_3$}
\label{tab:rel}
\end{table}

\begin{lemma}\label{lem:rel}
The elements $\delta$, $\varepsilon$ and $\Delta$ have the following relations.
\begin{enumerate}
\item
In $A(\arA_n)$, $\delta^{n+1} = \Delta^2 = \varepsilon^n$.

\item
In $A(\arB_n)$, $\delta^n = \Delta$.

\item
In $A(\arD_n)$,
$\delta^{n-1}=\Delta$ and
$\Delta^{\frac{2}{2\wedge n}}=\varepsilon^{\frac{n}{2\wedge n}}$,
hence $\delta^{2(n-1)} = \Delta^{2} = \varepsilon^{n}$.

\item
In $A(\arI_2(e))$,
$\delta^{\frac{e}{2\wedge e}}=\Delta^{\frac{2}{2\wedge e}}$,
hence $\delta^e=\Delta^2$.

\item
In $B(e,e,n)$,
$\delta^{n-1}=\Delta$ and
$\Delta^{\frac{e}{e\wedge n}}=\varepsilon^{\frac{n}{e\wedge n}}$,
hence $\delta^{e(n-1)} = \Delta^{e} = \varepsilon^{n}$.
\end{enumerate}
\end{lemma}

\begin{proof}
The relation in (i) for $A(\arA_n)$ is well-known,
and the relations
$\delta^n = \Delta$ in $A(\arB_n)$,
$\delta^{n-1}=\Delta$ in $A(\arD_n)$ and $B(e,e,n)$
and $\delta^{\frac{e}{2\wedge e}}=\Delta^{\frac{2}{2\wedge e}}$\
in $A(\arI_2(e))$ are immediate from the definitions of $\delta$ and $\Delta$.

We will prove only the relation
$\Delta^{\frac{e}{e\wedge n}}=\varepsilon^{\frac{n}{e\wedge n}}$ for $B(e,e,n)$.
Since the group $B(2,2,n)$ is the same as $A(\arD_n)$,
the relation $\Delta^{\frac{2}{2\wedge n}}=\varepsilon^{\frac{n}{2\wedge n}}$
in $A(\arD_n)$ is
a special case of the relation in $B(e,e,n)$ with $e=2$.

We will first prove the following identity:
\begin{equation}\label{eq:DeltaK}
\delta^k=(St_1St_2\cdots St_k)\cdot (s_{k+1}s_k\cdots s_3)\cdot t_{k+1},
\qquad k=1,2,\ldots,n-1,
\end{equation}
where the subscripts of $t$ are taken modulo $e$ as values between 1 and $e$.
When $k=1$, both sides of the equation (\ref{eq:DeltaK}) are identical.
Suppose the identity (\ref{eq:DeltaK}) holds for some $1\le k<n-1$.
We will show that it also holds for $k+1$.
Because $s_iS=Ss_{i+1}$ for $i=3,\ldots,n-1$, $t_1t_2=t_{k+1}t_{k+2}$, and
each $t_j$ commutes with $s_4,\ldots,s_n$, we have
\begin{align*}
(s_{k+1}\cdots s_3)\cdot t_{k+1}\cdot St_1t_2
&= (s_{k+1}\cdots s_3)\cdot t_{k+1}\cdot
  (s_n\cdots s_4)\cdot s_3\cdot t_{k+1} t_{k+2}\\
&= (s_{k+1}\cdots s_3)\cdot (s_n\cdots s_4)\cdot t_{k+1}
  \cdot s_3\cdot t_{k+1}\cdot t_{k+2}\\
&= (s_{k+1}\cdots s_3)\cdot (s_n\cdots s_4)\cdot s_3\cdot t_{k+1} \cdot s_3
\cdot t_{k+2}\\
&= (s_{k+1}\cdots s_3)\cdot S \cdot t_{k+1} \cdot s_3 \cdot t_{k+2}\\
&= S\cdot (s_{k+2}\cdots s_4) \cdot t_{k+1} \cdot s_3 \cdot t_{k+2}\\
&= S\cdot t_{k+1} \cdot (s_{k+2}\cdots s_4) \cdot s_3 \cdot t_{k+2}\\
&= S t_{k+1} \cdot (s_{k+2}\cdots s_4s_3) \cdot t_{k+2}.
\end{align*}
Therefore
\begin{align*}
\delta^{k+1}
&= \delta^k\cdot\delta
    =St_1St_2\cdots St_k \cdot (s_{k+1}s_k\cdots s_3) t_{k+1}\cdot St_1t_2\\
&= (St_1St_2\cdots St_{k+1})\cdot (s_{k+2}\cdots s_4s_3) \cdot t_{k+2}.
\end{align*}
This shows that the identity (\ref{eq:DeltaK}) holds for $k+1$.

\smallskip
When $k=n-1$, the identity (\ref{eq:DeltaK}) is the same as
$\delta^{n-1}=St_1St_2\cdots St_n$, hence we have
$$
\Delta=(St_1t_2)^{n-1}=\delta^{n-1}=St_1St_2\cdots St_n.
$$
Since the presentation of $B(e,e,n)$ is invariant under
the correspondence $t_j\mapsto t_{j+1}$, we have
$(St_jt_{j+1})^{n-1}=St_jSt_{j+1}\cdots St_{j+n-1}$.
Because $\delta=St_1t_2=St_jt_{j+1}$ for any $j\in\Z$,
we have
$$\Delta=St_jSt_{j+1}\cdots St_{j+n-1}
$$
for any $j\in\Z$.
Therefore for any $k\ge 1$
$$
\Delta^k=(St_1\cdots St_n)\cdot (St_{n+1}\cdots St_{2n})\cdots
(St_{(k-1)n+1}\cdots St_{kn}) = St_1\cdots St_{kn}.
$$
Because $\t_j=t_{j+e}$, we have for any $k\ge 1$
$$
\varepsilon^k=(St_1\cdots St_e)\cdot (St_{e+1}\cdots St_{2e})\cdots
(St_{(k-1)e+1}\cdots St_{ke})= St_1\cdots St_{k e}.
$$
By the above two identities, we have
$\Delta^{\frac{e}{e\wedge n}}=\varepsilon^{\frac{n}{e\wedge n}}$.
\end{proof}

The following theorem is an analogue for the groups
$A(\arB_n)$, $A(\arD_n)$, $A(\arI_2(e))$ and $B(e,e,n)$
of the Brouwer-Ker\'ekj\'art\'o-Eilenberg theorem for $A(\arA_n)$.

\begin{theorem}\label{thm:per-elt}
In the groups $A(\arA_n)$, $A(\arD_n)$ and $B(e,e,n)$,
every periodic element is conjugate to a power of $\delta$ or $\varepsilon$.

In the group  $A(\arB_n)$,
every periodic element is conjugate to a power of $\delta$.

In the group $A(\arI_2(e))$,
every periodic element is conjugate to a power of $\delta$ or $\Delta$.
\end{theorem}

\begin{proof}
We prove only the claim for $A(\arD_n)$ as
we can use the same argument  for the other groups.
As the center of $A(\arD_n)$ is a cyclic group generated by $\delta^{h'}$ with $h'$
a divisor of $h$, Theorem~\ref{thm:primitive} yields that
every periodic element in $A(\arD_n)$ is a power of a root of $\delta^{h'}$
and hence a power of a root of $\mu=\delta^{h}$.
Thus it suffices to show that every root of $\mu$ is conjugate to a power of
$\delta$ or $\varepsilon$.

Let $g$ be a $d$-th root of $\mu$ for a positive integer $d$.
By Theorem~\ref{thm:bes06a} and Lemma~\ref{lem:reg-num}, either $d|n$ or $d|2(n-1)$.
As $\delta^{2(n-1)}=\mu=\varepsilon^n$, there is a power of $\delta$ or $\varepsilon$
which is a $d$-th root of $\mu$.
By Theorem~\ref{thm:bes06a}, all the $d$-th roots of $\mu$ are conjugate to each other.
Therefore $g$ is conjugate to a power of $\delta$ or $\varepsilon$.
\end{proof}

\subsection{Precentrality}

This subsection investigates precentrality for periodic elements
in the classical and dual Garside structures on the groups
$A(\arA_n)$, $A(\arB_n)$, $A(\arD_n)$, $A(\arI_2(e))$ and $B(e,e,n)$.

Recall from Corollary~\ref{cor:Periodic_Tinf} that,
for a periodic element $g$ of a Garside group $G$,
$g$ is $p/q$-periodic if and only if
$\INF(g)=p/q$, $p\in\Z$, $q\in\Z_{\ge 1}$ and $p\wedge q=1$.

\begin{theorem}\label{thm:dual}
In the dual Garside structure on each of the groups
$A(\arA_n)$, $A(\arB_n)$, $A(\arD_n)$, $A(\arI_2(e))$ and $B(e,e,n)$,
every periodic element is either precentral or conjugate to a power of the Garside element $\delta$.
\end{theorem}

\begin{proof}
In the dual Garside structure on each group,
$\delta$ is the Garside element.
Hence the exponent $m$ of the minimal positive central power of
the Garside element $\delta$ is equal to $h'$ shown in Table~\ref{tab:rel}.

From Theorem~\ref{thm:per-elt} and Lemma~\ref{lem:precentral}, it suffices
to show that $\varepsilon$ is precentral in $A(\arA_n)$, $A(\arD_n)$ and $B(e,e,n)$,
and that $\Delta$ is precentral in $A(\arI_2(e))$.

In $A(\arA_n)$, $\varepsilon^n=\delta^{n+1}$ is the generator of the center,
hence $m=n+1$ and $\INF(\varepsilon)=(n+1)/n$.
Since $n+1$ and $n$ are coprime, $\varepsilon$ is $(n+1)/n$-periodic.
Therefore $\varepsilon$ is precentral in $A(\arA_n)$.

In $B(e,e,n)$, $\varepsilon^{n/(e\wedge n)}=\delta^{e(n-1)/(e\wedge n)}$
is the generator of the center,
hence $m=\frac{e(n-1)}{e\wedge n}$ and
$\INF(\varepsilon)=\frac{e(n-1)/(e\wedge n)}{n/(e\wedge n)}$.
As $e(n-1)\wedge n=e\wedge n$,
$\frac{e(n-1)}{e\wedge n}$ and $\frac{n}{e\wedge n}$ are coprime,
hence $\varepsilon$ is $\frac{e(n-1)}{e\wedge n}/\frac{n}{e\wedge n}$-periodic.
Therefore $\varepsilon$ is precentral in $B(e,e,n)$.

Since $A(\arD_n) = B(2,2,n)$, $\varepsilon$ is always precentral in $A(\arD_n)$.

In $A(\arI_2(e))$, $\Delta^{2/(e\wedge 2)}=\delta^{e/(e\wedge 2)}$
is the generator of the center, hence $m=\frac{e}{e\wedge 2}$ and
$\INF(\Delta)=\frac{e/(e\wedge 2)}{2/(e\wedge 2)}$.
As $\frac{e}{e\wedge 2}$ and $\frac{2}{e\wedge 2}$ are coprime,
$\Delta$ is $\frac{e}{e\wedge 2}/\frac{2}{e\wedge 2}$-periodic.
Therefore $\Delta$ is precentral in $A(\arI_2(e))$.
\end{proof}

In the classical Garside structure, every periodic element of $A(\arB_n)$ is precentral,
and every periodic element of $A(\arI_2(e))$ is either precentral or conjugate to a power of
the Garside element $\Delta$.
However it is not the case for the other groups.

\begin{theorem}\label{thm:artin}
In the classical Garside structure, the following hold.
\begin{enumerate}
\item
In $A(\arA_n)$, $\delta$ is precentral if and only if $n$ is even,
and $\varepsilon$ is precentral if and only if $n$ is odd.

\item
In $A(\arB_n)$, $\delta$ is always precentral.

\item
In $A(\arD_n)$, $\delta$ is precentral if and only if $n$ is even,
and $\varepsilon$ is always precentral.

\item
In $A(\arI_2(e))$, $\delta$ is always precentral.

\item
In $B(e,e,n)$, $\delta$ is precentral if and only if $n$ is a multiple of $e$,
and $\varepsilon$ is always precentral.
\end{enumerate}
\end{theorem}

\begin{proof}
In the classical Garside structure on each group,
$\Delta$ is the Garside element.
Let $m$ denote the exponent of the minimal positive central power of
the Garside element $\Delta$.

\medskip\noindent
(i)\ \
In $A(\arA_n)$, $\delta^{n+1}=\Delta^2=\varepsilon^n$ is the generator of the center,
hence $m=2$,
$$
\INF(\delta)= \frac{2}{n+1}\quad\mbox{and}\quad
\INF(\varepsilon)=\frac{2}{n}.
$$
Therefore, $\delta$ is precentral if and only if $n$ is even,
and $\varepsilon$ is precentral if and only if $n$ is odd.

\medskip\noindent
(ii)\ \
In $A(\arB_n)$, $\delta^n=\Delta$ is the generator of the center,
hence $m=1$ and $\delta$ is $1/n$-periodic.
Therefore $\delta$ is precentral.

\medskip\noindent
(iv)\ \
In $A(\arI_2(e))$, $\delta^{e/(e\wedge 2)}=\Delta^{2/(e\wedge 2)}$
is the generator of the center, hence  $m=\frac{2}{e\wedge 2}$.
Since $\frac{e}{e\wedge 2}$ and $\frac{2}{e\wedge 2}$ are coprime,
$\delta$ is $\frac{2}{e\wedge 2}/\frac{e}{e\wedge 2}$-periodic.
Therefore $\delta$ is precentral.

\medskip\noindent
(v)\ \
In $B(e,e,n)$, $\delta^{e(n-1)/(e\wedge n)} = \Delta^{e/(e\wedge n)}
= \varepsilon^{n/(e\wedge n)}$ is the generator of the center,
hence $m= \frac{e}{e\wedge n}$,
$$
\INF(\delta)= \frac{1}{n-1}
\quad\mbox{and}\quad
\INF(\varepsilon)=\frac{\ \frac{e}{e\wedge n}\ }{\ \frac{n}{e\wedge n}\ }.
$$
Since $\frac{e}{e\wedge n}$ and $\frac{n}{e\wedge n}$ are coprime,
$\varepsilon$ is always precentral.
And $\delta$ is precentral if and only if $m=\frac{e}{e\wedge n}=1$,
that is, if and only if $e|n$.

\medskip\noindent
(iii)\ \
Since $A(\arD_n)=B(2,2,n)$, $\delta$ is precentral if and only if $n$ is even,
and $\varepsilon$ is always precentral.
\end{proof}

\section{Discussions on some algorithmic problems}
This section discusses some algorithmic problems
concerning periodic elements in Garside groups.
As before, $G$ is a Garside group with Garside element $\Delta$,
and $\Delta^m$ is the minimal positive central power of $\Delta$.

\subsection{Periodicity decision problem}
\label{ssec:periodicity}

Let us consider the following problem.

\begin{quote}\em
Periodicity decision problem:
Given an element of a Garside group,
decide whether it is periodic or not.
\end{quote}

For Garside groups whose primitive periodic elements are well understood,
as for the groups $A(\arA_n)$, $A(\arB_n)$, $A(\arD_n)$, $A(\arI_2(e))$ and $B(e,e,n)$,
the periodicity decision problem is easy to solve.
From Theorem~\ref{thm:per-elt} and Table~\ref{tab:rel},
we can see the following.
\begin{enumerate}
\item
an element $g\in A(\arA_n)$ is periodic if and only if
either $g^n$ or $g^{n+1}$ is central;

\item
an element $g\in A(\arB_n)$ is periodic if and only if
$g^{n}$ is central;

\item
an element $g\in A(\arD_n)$ is periodic if and only if
either $g^{\frac{n}{2\wedge n}}$ or $g^{\frac{2(n-1)}{2\wedge n}}$ is central;

\item
an element $g\in A(\arI_2(e))$ is periodic if and only if
either $g^{\frac{2}{e\wedge 2}}$ or $g^{\frac{e}{e\wedge 2}}$ is central;

\item
an element $g\in B(e,e,n)$ is periodic if and only if
either $g^{\frac{n}{e\wedge n}}$ or $g^{\frac{e(n-1)}{e\wedge n}}$ is central.
\end{enumerate}
The centers of the above groups are cyclic generated by a power of the Garside element,
hence it is easy to decide whether a given element is central or not.

\smallskip
For arbitrary Garside groups, the periodicity decision problem
can be solved with a little more efforts.
From Lemma~\ref{lem:per-elt} and Proposition~\ref{prop:BasicTinf},
the following conditions are equivalent for an element $g$ of a Garside group $G$:
\begin{enumerate}
\item
$g$ is periodic;

\item
$g^q$ is conjugate to $\Delta^p$
for some $p,q\in\Z$ with $1\le q\le\Vert\Delta\Vert$;

\item
$g^{qm}=\Delta^{pm}$
for some $p,q\in\Z$ with $1\le q\le\Vert\Delta\Vert$.
\end{enumerate}
As the last two conditions can be checked by using the Garside structure,
the periodicity decision problem in Garside groups can be solved.

\subsection{Tabulation of primitive periodic elements}

For the groups $A(\arA_n)$, $A(\arB_n)$, $A(\arD_n)$,
$A(\arI_2(e))$ and $B(e,e,n)$,
primitive periodic elements were characterized in Theorem~\ref{thm:per-elt}.
However, for an arbitrary Garside group $G$,
we know by Theorem~\ref{thm:primitive} only that
every primitive periodic element in $G$ is a $k$-th root of $\Delta^m$ for
some $1\le |k|\le m \|\Delta\|$.
On the other hand, periodic elements have summit canonical length 0 or 1.
Hence every primitive periodic element is conjugate to an element of the form
$\Delta^u a$ for $-m\le u\le m$ and $a\in[1,\Delta)$.
Therefore there are only finitely many primitive periodic elements in $G$ up to conjugacy.
We consider the following problem.

\begin{quote}\em
Primitive periodic element tabulation:
Given a Garside group $G$,
make a list of primitive periodic elements
such that each primitive periodic element of $G$
is conjugate to either exactly one element in the list or its inverse.
\end{quote}

We solve the above problem in Proposition~\ref{prop:root}
by using the following solution to the root problem in Garside groups.

\begin{theorem}[\cite{Sty78,Sib02,Lee07}]\label{thm:root}
Let $G$ be a Garside group.
There is a finite-time algorithm that,
given an element $g\in G$ and an integer $k\ge 1$,
decides whether there exists $h\in G$ with $h^k=g$,
and then finds such an element $h$ if one exists.
\end{theorem}

The above theorem was proved for braid groups by Sty\v shnev~\cite{Sty78},
and for Garside groups by Sibert~\cite{Sib02} and Lee~\cite{Lee07}.

\begin{proposition}\label{prop:root}
Given a Garside group $G$, there exists a finite-time algorithm
that makes a list of primitive periodic elements
such that each primitive periodic element of $G$
is conjugate to either exactly one element in the list or its inverse.
\end{proposition}

\begin{proof}[Sketchy proof]
The algorithm performs the sequential steps below.
\begin{enumerate}
\item[Step 1.]
Compute all roots of $\Delta^m$ (up to inverse and conjugacy):
for each element $h$ of the form $\Delta^u a$ for $0\le u\le m$ and $a\in[1,\Delta)$,
decide whether $h^k=\Delta^m$ for some $1\le k\le m\Vert\Delta\Vert$.

\item[Step 2.]
Let $H=\{h_1,\ldots,h_N\}$ be the set of all roots of $\Delta^m$
obtained from the above step.
As the conjugacy problem is solvable in Garside groups,
we can partition the set $H$ into conjugacy classes,
and then select one element from each conjugacy class.
In this way, we obtain a subset $H'$ of $H$ such that
each root of $\Delta^m$ is conjugate to either exactly one element of $H'$ or its inverse.

\item[Step 3.]
For each element $h$ of $H'$, decide whether $h$ has
a $k$-th root for $2\le k\le m\Vert\Delta\Vert$,
and remove $h$ from $H'$ if it does.
(By Theorem~\ref{thm:root}, this can be done in a finite number of steps.)
Let $H''$ be the resulting set.
Then each primitive periodic element in $G$ is conjugate to either exactly one element of $H''$
or its inverse.
\end{enumerate}
\end{proof}

\subsection{Conjugacy problem for periodic elements}\label{ssec:CDP:CSP}

Observe that
the relations in the presentations of Artin groups and the braid group
of complex reflection groups are all homogeneous.
Therefore, the exponent sum of an element, written as a word in
the generators and their inverses, is well defined.
The exponent sum is invariant under conjugacy.

Let us consider the CDP and CSP for periodic elements in Garside groups.

First, we shall see that the exponent sum is a complete invariant
for the conjugacy classes of periodic elements in the groups
$A(\arA_n)$, $A(\arB_n)$, $A(\arD_n)$, $A(\arI_2(e))$ and $B(e,e,n)$,
hence the CDP is easy in those groups.
To establish this, we need the fact that the roots of periodic elements are
unique up to conjugacy in these groups.
This was proved by Bessis, see Theorem~\ref{thm:bes06a}.
Because his theorem is stated only for the roots of $\Delta^m$,
we prove the following lemma for completeness of the paper.

\begin{lemma}\label{lem:RootPer}
Let $G$ be a Garside group such that, for any $k\ge 1$,
the $k$-th root of $\Delta^m$, if exists, is unique up to conjugacy.
Then, for any periodic elements $g_1, g_2 \in G$ and for any nonzero integer $k$,
$g_1^k=g_2^k$ implies that $g_1$ is conjugate to $g_2$.
\end{lemma}

\begin{proof}
Choose any periodic elements $g_1, g_2 \in G$ and any nonzero integer $k$.
Suppose $g_1^k=g_2^k$.
Let $g_1$ be $p/q$-periodic, then so is $g_2$ because $\INF(g_1)=\INF(g_2)$.
If $p=0$, there is nothing to do. Let $p\ne 0$.
Applying Lemma~\ref{lem:PerCyclic} to $g_1$ and $g_2$, we have the following:
there are $r, s\in\Z$ with $pr+qms=p\wedge m$;
let $h_i = g_i^r \Delta^{ms}$ for $i=1,2$;
then $h_1$ and $h_2$ are $\frac{qm}{p\wedge m}$-th roots of $\Delta^m$.
Hence $h_1$ is conjugate to $h_2$ by the hypothesis on $G$.
On the other hand, $g_1 = h_1^{\frac{p}{p\wedge m}}$ and $g_2 = h_2^{\frac{p}{p\wedge m}}$
by Lemma~\ref{lem:PerCyclic}.
Therefore $g_1$ is conjugate to $g_2$.
\end{proof}

\begin{proposition}\label{prop:ExpSumConj}
Let $G$ be one of the Garside groups $A(\arA_n)$,
$A(\arB_n)$, $A(\arD_n)$, $A(\arI_2(e))$ and $B(e,e,n)$.
Let $g_1$ and $g_2$ be periodic elements in $G$.
Then, $g_1$ and $g_2$ are conjugate if and only if
they have the same exponent sum.
\end{proposition}

\begin{proof}
Suppose that $g_1$ and $g_2$ have the same exponent sum.
There is an integer $k\ge 1$ such that
both $g_1^k$ and $g_2^k$ belong to $\langle\Delta^m\rangle$.
As $g_1$ and $g_2$ have the same exponent sum,
we have $g_1^k = g_2^k$.
For any $d\ge 1$, the $d$-th root of $\Delta^m$,
if exists, is unique up to conjugacy by Theorem~\ref{thm:bes06a}.
Therefore $g_1$ and $g_2$ are conjugate by Lemma~\ref{lem:RootPer}.
The converse is obvious.
\end{proof}

From the above proposition, it is easy to solve the conjugacy decision problem
for periodic elements in the groups
$A(\arA_n)$, $A(\arB_n)$, $A(\arD_n)$, $A(\arI_2(e))$ and $B(e,e,n)$.

\smallskip
Now we will consider the conjugacy search problem for periodic elements in those groups,
using the dual Garside structure.

In the groups $A(\arB_n)$ and $A(\arI_2(e))$,
the CSP for periodic elements is easy to solve.
In $A(\arB_n)$, every periodic element is conjugate
to $\delta^k$ for some $k\in\Z$.
Hence the super summit set is of the form $\{\delta^k\}$
since $\delta$ is the Garside element.
As for $A(\arI_2(e))$,
the dual presentation is
$$
A(\arI_2(e))=\langle t_1,\ldots,t_e\mid
t_1 t_2 = t_2 t_3 =\cdots = t_{e-1}t_{e}= t_{e}t_{1}\rangle.
$$
Since $\delta=t_1t_2$ is the Garside element,
the set of simple elements is $[1,\delta]=\{1,t_1,t_2,\ldots,t_e,\delta\}$,
hence the super summit set of a periodic element
is of the form either $\{\delta^k\}$
or a subset of $\{\delta^k t_i:i=1,\ldots,e\}$.
Therefore for both groups $A(\arB_n)$ and $A(\arI_2(e))$,
the super summit set of a periodic element
is very small, hence the CSP is easy to solve.

In the groups $A(\arA_n)$, $A(\arD_n)$ and $B(e,e,n)$,
every periodic element is conjugate to a power of $\delta$ or $\varepsilon$.
Since $\delta$ is the Garside element in the dual Garside structure,
the CSP is easy to solve for conjugates of $\delta^k$, $k\in\Z$,
because their super summit set consists of a single element.
Therefore it is enough to consider the conjugates of powers of $\varepsilon$.
Given a conjugate $\alpha$ of $\varepsilon^k$ for a nonzero integer $k$,
it is easy to compute a nonzero integer $r$,
by Lemmas~\ref{lem:BCMW-a} and~\ref{lem:BCMW}, such that
$\varepsilon^{kr}$ is slim and the CSP for $(\alpha, \varepsilon^k)$ is equivalent to
the one for $(\alpha^r, \varepsilon^{kr})$.
By Theorem~\ref{thm:dual} and Lemma~\ref{lem:precentral},
every power of $\varepsilon$ is precentral.
Therefore we may assume that
the given periodic elements are slim and precentral
so that we can use Theorem~\ref{thm:pa-cy} that the super summit
set is closed under any partial cycling.

For periodic elements in arbitrary Garside groups,
we do not know whether the CDP is easier than the CSP.
The CSP for periodic elements looks easier than for arbitrary elements,
because the super summit elements are of
the form $\Delta^k a$ for $a\in[1,\Delta)$,
and because in case the periodic elements are precentral
the super summit set can be assumed to be closed under any partial cycling.

\subsection*{Acknowledgements}

The authors are grateful to the referee
for valuable comments and suggestions.
The first author was supported by Basic Science Research Program
through the National
Research Foundation of Korea (NRF) funded by the Ministry of Education,
Science and Technology (2009-0063965).
The second author was supported by Konkuk University
and by Basic Science Research Program through the National
Research Foundation of Korea (NRF) funded by the Ministry of Education,
Science and Technology (2010-0010860).

\end{document}